\documentclass{elsarticle}

\usepackage[utf8]{inputenc}
\usepackage[T1]{fontenc}

\usepackage{anysize}

\usepackage{cancel}
\usepackage{amssymb}
\usepackage{array,longtable}
\usepackage{amsmath}
\usepackage{amsthm}
\usepackage{xcolor}
\usepackage{graphicx}
\usepackage{rotating}
\usepackage{algorithm}
\usepackage{multirow}
\usepackage{multicol}
\usepackage{algorithmic}
\usepackage{graphicx}
\usepackage{enumerate}
\usepackage{enumitem}
\usepackage{url}
\usepackage{hhline}
\usepackage{subcaption}
\usepackage{mathtools}
\usepackage{rotating}
\usepackage{empheq}
\usepackage{hyperref}
\hypersetup{colorlinks=true,filecolor=blue,linkcolor=blue,citecolor=blue,urlcolor=blue,pdfnewwindow=true}
\newcolumntype{H}{>{\setbox0=\hbox\bgroup}c<{\egroup}@{}}

\usepackage{natbib}
 \bibpunct[, ]{(}{)}{,}{a}{}{,}

\usepackage{algorithm}
\usepackage{algorithmic}
\usepackage{textcomp}
\usepackage{longtable}

\biboptions{authoryear}

\newcolumntype{d}[1]{D{.}{.}{#1}}

\newcommand{\BB}{{\cal B}}
\newcommand{\C}{{\cal C}}

\renewcommand{\H}{{\cal H}}
\newcommand{\I}{{\cal I}}
\newcommand{\J}{{\cal J}}\newcommand{\K}{{\cal K}}
\renewcommand{\L}{{\cal L}}
\newcommand{\M}{{\cal M}}
\newcommand{\N}{{\cal N}}

\renewcommand{\P}{{\cal P}}
\newcommand{\Q}{{\cal Q}}

\newcommand{\bea}{\begin{eqnarray}}
\newcommand{\eea}{\end{eqnarray}}

\journal{EURO Journal on Transportation and Logistics}

 \newcommand{\bcb}{\color{black}}

\begin{document}

\begin{frontmatter}

\title{Cross-Dock Door Design under Uncertainty: A two-stage DRO-based lower- and upper-bounding scheme}

\author[URJC]{Laureano F. Escudero}
\ead{laureano.escudero@urjc.es}
\author[UPV-EHU]{M. Araceli Gar\'in\corref{cor}}
\ead{mariaaraceli.garin@ehu.eus}
\author[UPV-EHU1]{Aitziber Unzueta}
\ead{aitziber.unzueta@ehu.eus}
\cortext[cor]{Corresponding author}

\address[URJC]{Area of Statistics and Operations Research, Universidad Rey Juan Carlos, URJC, M\'ostoles (Madrid), Spain.}
\address[UPV-EHU]{Quantitative Methods department, Universidad del Pa\'{\i}s Vasco, UPV/EHU, Bilbao (Bizkaia), Spain.}
\address[UPV-EHU1]{Applied Mathematics department, Universidad del Pa\'{\i}s Vasco, UPV/EHU, Bilbao (Bizkaia), Spain.}

\begin{abstract}
The stochastic cross-dock door design problem entails determining the
number of doors and their nominal capacities under uncertainty.
The inbound flow of commodities from origin nodes is assigned to the entry
doors consolidated in the platform, and the outbound flow is
assigned to the exit doors to be delivered to the destination nodes.
This problem combines three high computational difficulties, namely, NP-hard quadratic combinatorics, uncertainty in
the main parameters, and ambiguity in their probability
distribution. Distributionally robust optimization is considered to
deal with these uncertainties. A  two-stage mixed binary quadratic
model is presented, where the first stage decisions are
related to the design of the platform and the second stage ones are related
to the assignment of the commodity flow to the doors in the members of the ambiguity set.
The goal is to minimize
the highest total cost in the ambiguity set, subject to the constraint
system for each of those members.
In addition to the risk-neutral version, a risk-averse formulation is presented.
Given the difficulty of this problem,  a
\textit{min-max} matheuristic scheme based on a scenario cluster
decomposition is proposed for obtaining lower and upper bounds.
A computational study is conducted to compare the solutions provided by
the straightforward use of the state-of-the-art solvers CPLEX and
Gurobi, as well as to validate the proposed matheuristic scheme.
\end{abstract}

\begin{keyword}
Combinatorial optimization; cross-dock door design; distributionally
robust optimization; risk-neutral; risk-averse.
\end{keyword}

\end{frontmatter}
\section{Introduction and motivation}\label{sec:intro}
\begin{figure}[h]
\begin{center}
\includegraphics[width=12cm]{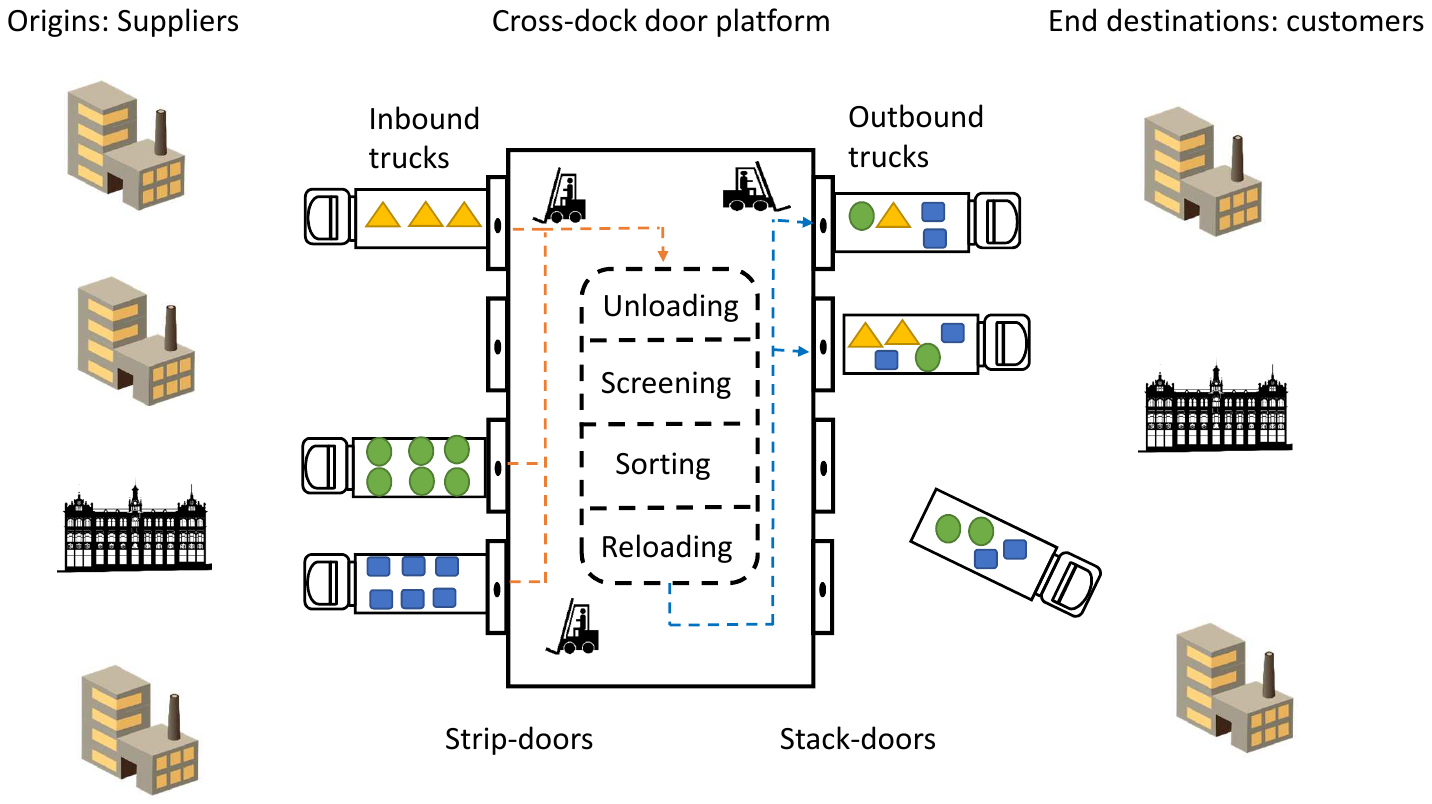}
\vskip -1.7cm
\caption{CD platform}\label{fig:cross-dock}
\end{center}
\end{figure}

Given a network in which commodities must be delivered from origin nodes to destination nodes, the cross-docking problem consists of using a facility that serves as a consolidation point to streamline the distribution within the supply chain.
The origin nodes can deliver the commodities to the cross-dock
platform, through entry points referred to as {\it strip} doors, so that after
being classified by type and destination they exit (usually
in smaller quantities) to the destination nodes, through  exit doors,
known as {\it stack} doors.
 The classical operation of the platform is depicted in Fig. \ref{fig:cross-dock}, taken from \cite{cdap24}.

One of the optimization problems in cross-docking is the so-called
Cross-Dock Door Assignment Problem (CDAP). The objective of this
problem is to assign the commodity flows to the strip and stack doors in a way that minimizes operational costs.
  It is a deterministic problem,
so, the set of available doors is given as well as their respective capacities.
Even considering these sets and parameters as deterministic ones, CDAP is a difficult binary quadratic combinatorial problem.
Its constraint system is composed of two {\it Generalized Assignment Problems}:
one is related to the strip doors and the other is related to the stack doors,
to which the origin nodes and the destination nodes are to be assigned, respectively.
Those two problems are linked by quadratic terms in the objective function, see \cite{Guignard12}.

Another important and also difficult combinatorial problem in
cross-docking is the so-called Cross-Dock Door Design Problem (CDDP),
the subject of this paper. The goal of this stochastic problem is
twofold: first, to select the doors and their nominal
capacities (strategic decisions); and second, to assign commodity flows
to the strip and stack doors (operational decisions). Here, the set of available strip and stack doors is not given, nor the capacities of the doors. On
the other hand, the sets of origin and destination nodes
are uncertain and, so, the volumes of the commodity flow from the origins to
the destinations is also uncertain. The other sources of uncertainty lie
in the handling and transportation costs through the platform as well
as the door capacity disruption. There are two types of node--door
assignments: the standard procedure and the {\it  outsourcing} one.
The latter must be considered when either no strip door is available for an inbound commodity flow or
no stack door is available for an outbound commodity flow.
Outsourcing is required in both cases. Notice that CDDP could
be formulated as a two-stage stochastic model,
where the assignment is performed under a unique finite set of scenarios,
which represent the realizations of the uncertain parameters,  see \cite{cddp-ts24}.  In particular, the operational decisions under each second-stage
scenario of CDDP correspond to the solution of a CDAP.

\bigskip
\noindent\textbf{CDDP-DRO description as a subject of this work}

While the stochastic approach to the problem is very attractive when
the scenario set closely approximates the true underlying probability
distribution, the alternative proposed in this paper addresses the case
in which the true distribution is unknown, which adds an additional difficulty to deal with.

The two-stage stochastic approach for CDDP introduced in this paper considers
{\it distributionally robust optimization} (CDDP-DRO), see \cite{Khun25}, to get
a better accuracy in dealing with the uncertainty.
 Rather than using DRO, this work presents a two-stage  model based on
the Wasserstein metric to capture  uncertainties, that themselves are grouped in the so-called ambiguity sets.

On the other hand, the risk-neutral formulation underestimates the costs by averaging the extreme
scenarios, commonly referred to as black swans.  So, a high negative
impact in the objective function cost should be prevented for those
black swan scenarios that have a low probability and a high cost.
 In this context, when the platform lacks available strip or stack
doors to accommodate inbound or outbound commodity flows, it becomes
necessary to outsource the handling operations. However, due to the
high associated costs, such outsourcing should be avoided whenever possible.
 For that purpose, an outsourcing risk-averse  approach is introduced, in addition to the risk-neutral model,
whose aim is to prevent solutions that require outsourcing.

All of the above elements considered contribute to the overall complexity of the problem, thereby requiring a decomposition matheuristic that  exploits the specific structure of the
model. The algorithm that we propose introduces a \textit{min-max} strategy for obtaining a (hopefully, good) robust feasible solution and
a \textit{max-max} strategy for obtaining a (hopefully, strong) lower bound on the optimal solution.
Those strategies are based on the decomposition of each
member of the ambiguity set into scenario clusters.
As far as we know, the challenging problem that results has not been addressed before, although it has a very broad application field.

\bigskip
\noindent\textbf{Literature review}

The treatment of the cross-docking problem is a young discipline from
an academic point of view: as a matter of fact, most of the literature
on this type of problem has been published over the last 25 years.
For good overviews of applications of cross-docking, see \cite{Akkerman22}, among others.
For different approaches to deterministic CDAP solving, see \cite{cdap24,Gelareh20,Guignard20,Monemi23} and \cite{Nassief18},
mainly for matheuristic algorithms, and an extensive review in \cite{Goodarzi20}.
Comprehensive reviews of CDAP under uncertainty have been recently published
in \cite{ArdakaniFei20} and \cite{BuakumWisittipanich19}.

Two stochastic quadratic models are presented in \cite{Soanpet12} for CDDP,
where the  uncertainty lies in the centres' capacity disruption that is assumed to follow the Normal distribution.
The risk-averse approach is the Chance-Constrained Programming (CCP), also so-called first-order stochastic dominance measure.
A binary linear optimization model is presented in \cite{Essghaier21} for assigning trucks to doors in a given platform,
where a collaboration of supply competitors is allowed in a fuzzy CCP-based approach.
\cite{Goodarzi20} present an approach for capacity disruption issues
in the truck scheduling problem.
A column-and-constraint generation exact algorithm is proposed.
See \cite{Xi20} for uncertainty in the arrival and operational times.
\cite{Gallo22} present a two-stage stochastic Mixed Integer Linear Programming (MILP) model for scheduling trucks to doors.
A genetic matheuristic algorithm is considered.
An efficient heuristic, so-called the Stability Approach, is presented in \cite{Fonseca24} for solving multi-dock truck scheduling problems,
where the uncertainty lies in the truck release dates.

It is worth pointing out that, except for \cite{cddp-ts24}, the
reviewed  literature on both CDAP under uncertainty and  CDDP does not
consider the specific capacity of each door in the identification and assignment of individual doors.
In fact, only the total capacity of each platform is considered, which implies an approximation in the solution.

An important piece of information to deal with in practical
applications as CDDP is the probability distributions followed by the uncertainty in the main parameters.
It is usually unknown, being the motivation for DRO,
where a set of possible distributions is dealt with as a counterpart of the true underlying one;
\cite{Khun25} present a good survey on DRO key elements in a unified manner.
\cite{LiLiu25} editorialize about some models and methods for generating the ambiguity set.
To the best of our knowledge, \cite{Scarf57} is the first to consider DRO (although without using this name)
 for solving the famous news vendor problem, where the mean and variance are assumed to be known.
There is a broad literature where given bounds are considered for the
distance of possible distributions from an empirical one.
The variety of measures to define that distance includes
the popular Wasserstein distance introduced in \cite{Kantorovich42}, also known as the translation of masses or Earth mover model;
the nested distance (\cite{PflugPichler14});
the $\phi$ divergences as versions of the $\chi^2$ distance, the Kullback--Leibler divergence, the Hellinger distance, and the Burg entropy
(see \cite{Philpott18} and \cite{Bayraksan24}, among others);
the $\ell_1$ and $\ell_{\infty}$ norms defining distances (\cite{Xie20,Wang23});
$\phi$ divergence functions introduced in \cite{CaunhyeAlem23}, which are based on infima convolution, piecewise linearity, and the Moreau--Yosida regularization.

The Wasserstein metric has been widely employed in applications of DRO; see, for example, \cite{Beesten24,Duque22,ElTonbary24,ErdoganIyenhar06,GaoKleywegt22,Jackiewicz23,ShenJiang23,Yuan25}. Several extensions have also been studied, including the barycenter Wasserstein distance (\cite{MohajeringKhun18,Arrigo22}), the $\infty$-Wasserstein distance (\cite{Bertsimas22}), the $1$-norm Wasserstein distance (\cite{Agra25,Byeon25}), and the Sinkhorn-regularized Wasserstein distance (\cite{Seyedi25,Ou26,Wang26}), among others.

Appendix A in the Supplementary file provides an overview of DRO applications.
It is worth noting that the DRO approach introduced in this work is new to the cross-docking platform design problem.

The literature on {\it risk-averse} DRO is very scarce and, as far as we know, it is basically restricted to Conditional Value-at-Risk (CVaR) and CCP.
\cite{Jackiewicz23} present a robust approach with uncertainty in the coefficients of the objective function.
CVaR is proposed as the risk-averse measure in the model,
where the Wasserstein metric defines the ambiguity set in the neighbourhood ball centred on an empirical distribution.
A risk-averse DRO approach is presented in \cite{Lin22} for principal--agent contract design.
The model considers the worst case where the objective function also includes the CVaR risk-averse measure.
\cite{Gao24} present a two-stage MILP DRO scheme, where the risk aversion is dealt with via a CCP approach.
Another DRO CCP is proposed in \cite{Sarmadi25},
for the selecting the best locations of rescue facilities for an uncertain occurrence of  a disaster.
 See \cite{Zhang25} for another DRO CCP approach.
\cite{Xin25} present a DRO approach for dealing with the location and
number of facilities for the disposal  of construction and demolition waste.
A two-stage MILP is proposed where a hybrid risk-averse measure composed of CVaR and CCP is considered.
\cite{Ou26} present a continuous nonlinear DRO CCP approach where the
Sinkhorn-based regularized Wasserstein distance is considered.
A DRO CVaR approach is dealt with in \cite{Zhai25} for the classical newsvendor problem.
For other DRO-CVaR approaches, see \cite{Ghaoui03,Huong08,Tong09} and \cite{Zhu09}.

The risk-averse measure so-called second-order stochastic dominance was
introduced in \cite{dentchevaruszczynski03}, see also \cite{Eppen89,Shapiro09,Shapiro21}, among others.
There are few works that consider MILP for that type of measure, see \cite{Gollmer11,ejor16}, \cite{ejor19,ejor20,Michelli26}.
DRO with second-order stochastic dominance constraints was introduced in \cite{dentchevaruszczynski10}, establishing its optimality conditions and considering continuous variables for the one-stage setting.
 See also \cite{Guo17,PengDelage20,Shapiro21,Zhang21,Hosseini-Nodeh22,Mei22}, where the Wasserstein distance is considered in most of the works.
In any case, to the best of our knowledge, none of the DRO approaches in the literature addresses the specific structural features of the cross-docking model.

Finally, regarding solution procedures for DRO, moderately sized problems can be solved using off-the-shelf solvers. For large instances, Benders-type decomposition algorithms are commonly employed, sometimes with extensions. See \cite{Jin26} for a specific decomposition approach for problems with a particular structure.

\bigskip
\noindent\textbf{Our main contributions and organization of this work}

The main contributions of this work are as follows:
\begin{enumerate}
\item A {\it risk-neutral two-stage CDDP-DRO model} is proposed,
    where quadratic expressions are required for modeling the cost interaction of the strip and stack doors in the problem.
    The model explicitly deals with the ambiguity set to cope with the uncertainty in the objective function and the constraint system.
    The goal consists of minimizing the highest scenario  expected cost in the ambiguity set.
\item An \textit{outsourcing risk-averse model} that incorporates a \textit{second-order stochastic dominance} constraint system is proposed for the two-stage CDDP-DRO mixed integer quadratic model
as a counterpart of the risk-neutral one.
A set of policy profiles can be considered for the modeler-defined key thresholds in the CDDP cost function.
The impact of each cost policy on the scenario highest expected cost in the ambiguity set is modeled.    For each policy, the bi-objective function to minimize consists of the weighted highest cost in the ambiguity set and the penalization of the scenario expected cost surplus over the  threshold related to the ambiguity set member that provides the highest cost.		In this way, the DRO solution is protected against the surplus over the cost thresholds all along the given policy profile.

\item Given the special structure of those highly combinatorial versions of CDDP and the high dimensionality of realistic instances, a mathematically equivalent model is introduced; it is so-called the {\it Linearized mixed Integer Programming} model.

\item A \textit{min-max} {\it matheuristic} approach is presented.
  It is based on a sophisticated {\it scenario cluster decomposition} that exploits the special structure of the problem.
  As a result, it obtains a set of candidate first stage feasible solutions as well as a lower bound on the optimal  one of the original model.
  A second stage solution is also obtained for each candidate first stage one,
  so that an upper bound is provided as the incumbent CDDP-DRO solution.
  \end{enumerate}

A computational experiment is reported to compare the solutions provided by
the straightforward use of the state-of-the-art commercial solvers CPLEX and Gurobi, as well as to validate the proposed matheuristic scheme.
To enable this experimentation, the ambiguity set of the uncertain parameters has been constructed using the Wasserstein metric.
Given a modeler-defined set of candidate distributions, the scheme derives the so-called Wasserstein ball from
the projections of the appropriate perturbations onto the distribution
functions of the second-stage realizations.
Those projections should be within a given radius in the Wasserstein distance from the empirical distribution.

\bigskip
The remainder of this paper is organized as follows. Section
\ref{sec:math} details the proposed model's specific robust formulation. Section
\ref{sec:RN} introduces the risk-neutral version of the two-stage CDDP-DRO model and presents the matheuristic scheme to obtain lower and upper bounds.
Section \ref{sec:ORA} presents the outsourcing-risk aversion functional, along with the additional details required by the proposed decomposition scheme.
Section \ref{sec:results} presents the numerical results and the empirical findings.
Section \ref{sec:con} draws some conclusions and outlines a future research agenda for generalizing the decomposition algorithm for solving any two-stage MILP problem under uncertainty.
Appendices B and C in the Supplementary file introduce the procedure for generating the ambiguity set, as well as some details of the DRO solutions obtained in the computational experiment.

\section{Mathematical formulation}\label{sec:math}
The aim of the two-stage mixed integer quadratic DRO models presented in this work for CDDP is to
obtain an incumbent feasible solution of the cross-docking platform design problem for alternatively satisfying the risk-neutral and the outsourcing risk-averse measures, providing also a quality guarantee of the solution.
So, the need for a robust treatment of the uncertainty adds difficulty to the computational treatment of the stochastic combinatorial model
--which is already NP-hard with a quadratic objective function --.

\subsection{Risk-Neutral two-stage stochastic metamodel}
 The risk-neutral metamodel can be expressed as
  \begin{eqnarray}\label{1bs_RN}
 &\min_{\delta} c^{T}\delta + \mathbb{E}_{\mathbb{P}}[\Q(\delta,
  \tilde{\xi})] &\mbox{s.to} ~~ \delta \in \Delta
  \end{eqnarray}
        where $c\in \mathbb{R}^{z}$ and $\delta\in \{0,1\}^{z}$ are the vectors of the first-stage objective function cost parameters and variables, respectively, $\Delta$ gives the feasible region, and $\mathbb{E}_{\mathbb{P}}[\Q(\delta,\tilde{\xi)}]$ is the expected value of the second-stage recourse function in the uncertain factors $\tilde{\xi}$.
The $\delta$-variables represent the decisions to be made on the strip and stack doors selection as well as the related capacity,
before the uncertainty in $\tilde{\xi}$ is unveiled.
For a given realization $\xi$ of $\tilde{\xi}$, the second-stage binary quadratic function $\Q(\delta,\xi)$ can be expressed as
\begin{subequations}\label{2bs_RN}
\bea
\Q(\delta, \xi)&=& \min_\nu  \nu_1 (\xi)^{T} G(\xi) \nu_2(\xi)\\
&&\mbox{s.to}~~ A(\xi)\delta + B(\xi)\nu(\xi) \le b(\xi)\\
&& \nu(\xi)=(\nu_1(\xi), \nu_2(\xi)), \,\, \nu_1(\xi) \in \{0,1\}^m, \,\, \nu_2(\xi) \in \{0,1\}^n,
\eea
\end{subequations}
where the decisions on the inbound and outbound commodity flows assigned to the doors are represented by the vectors of variables $\nu_1(\xi)$ and $\nu_2(\xi)$, respectively, $G(\xi)$ is the matrix of the cost interaction of the duplex given by [(inbound flow, strip doors), (outbound flow, stack doors)],
$A(\xi)$ and $B(\xi)$ are the constraint matrices and $b(\xi)$ is the right-hand-side vector, all with the appropriate dimensions.
The random vector $\tilde{\xi}$ is defined on a probability space and $\mathbb{P}$ represents its associated probability measure.

\subsection{Modelling uncertainty based on Wasserstein metric}
Let $\P(\Xi)$ denote the set of all probability distributions defined on $\Xi$ for the uncertain factors  $\tilde{\xi}$, and
$\P \subset \P(\Xi)$ is termed as an {\it ambiguity set} comprising all associated probability distributions that are compatible with the available information.
That information, so-called {\it nominal} or (reference) distribution $\mathbb{P}_{N}$ represents an empirical distribution that is not  guaranteed to be  the true probability one. It is worth noting that the ambiguity set can be constructed as a ball within the space of probability distributions
by exploiting a probability distance function. Note that, via adjusting the
 radius of the Wasserstein ball, the degree of conservatism of the
 distribution uncertainty can thus be controlled by the decision
 makers. Additionally, the worst-case expectation over the Wasserstein
 ball can be solved efficiently using off-the-shelf convex
 optimization techniques. So, the Wasserstein
 ambiguity set $\P$ is defined as
        \begin{eqnarray}\label{Wam-set}
 \P=\{\mathbb{P}: d_{W^{r}}( \mathbb{P}, \mathbb{P}_{N})\le \theta \}
 \end{eqnarray}
        where $\theta>0$ is  the Wasserstein
 radius, and $d_{W^{r}}(,)$ is the $r-$Wasserstein
 metric, defined as
 \begin{eqnarray}\label{Wass-metric}
 d_{W^{r}}( \mathbb{P}_{1}, \mathbb{P}_{2})=\inf_{\Pi}\left(
   \int_{\Xi x \Xi}||\xi_{1}-\xi_{2}||^{r}
   ~ \Pi(d\xi_{1},d\xi_{2})\right )^{1\over r}
\end{eqnarray}
 such that $\Pi$ is the joint
   distribution of $\tilde{\xi}_{1}$ and   $\tilde{\xi}_{2}$ with
   marginals $\mathbb{P}_{1}$ and $ \mathbb{P}_{2}$, respectively.

\subsection{Robust model formulation}\label{dro-form}
        A type of DRO objective function, as an alternative to \eqref{2bs_RN}, can be expressed such that
the first stage cost is not considered in the highest expected one in the ambiguity set members, see \cite{Beesten24}.
    However, we argue  that the first stage cost is better protected by
 including it in the robustness, so,  the
distributionally robust two-stage stochastic model is formulated as
  \begin{eqnarray}\label{DRO_RN}
 &\min_{\delta}\sup_{\mathbb{P}\in \P}\{ c^{T}\delta + \mathbb{E}_{\mathbb{P}}[\Q(\delta,
  \tilde{\xi})]\} &\mbox{s.to} ~~ \delta \in \Delta
  \end{eqnarray}

To formulate the risk-neutral two-stage DRO metamodel so that \eqref{DRO_RN} can be implemented, we introduce the auxiliary continuous variable $u$, representing the worst-case (i.e., highest)  scenario expected cost in the ambiguity set $\P$.
The uncertainty realizations are discretized into a set of scenarios
 ($\Omega_{\mathbb{P}}$) for each member of the ambiguity set,
 $\mathbb{P}\in \P$.  The so-called DRO-RN metamodel can then be expressed as
\begin{subequations}\label{robust_RN}
 \begin{eqnarray}
   &&\min u\\
  \mbox{s.to}& & \label{DRO_RN-b} c^{T}\delta + \sum_{\omega\in\Omega_{\mathbb{P}}}w^{\omega}\nu_1^{\omega T} G^{\omega} \nu_2^{\omega} \le u
                                                                               ~~\forall \mathbb{P}\in \P\\
   && \label{DRO_RN-c} A^{\omega}\delta + B^\omega \nu^{\omega} \le b^{\omega} ~~ \forall \omega\in\Omega_{\mathbb{P}}, ~\mathbb{P}\in \P\\
   && \label{DRO_RN-d}\delta \in \{0,1\}^{z}, \,\,
   \nu^\omega=(\nu_1^\omega, \nu_2^\omega), \,\, \nu_1^\omega \in \{0,1\}^m, \,\, \nu_2^\omega \in \{0,1\}^n
                                                                               ~~ \forall \omega\in\Omega_{\mathbb{P}}, ~~\mathbb{P}\in\P,
\end{eqnarray}
\end{subequations}
where $w^\omega$ is the weight of scenario $\omega$ in $\Omega_{\mathbb{P}}$, such that $\sum_{\omega_{\mathbb{P}}}w^\omega=1$.

Formulation \eqref{robust_RN} is risk-neutral in the sense that it is
only concerned with the expected recourse objective function.
 As discussed in the introduction, the DRO literature scarcely addresses the use of risk-averse  measures such as CVaR, CCP and {\it stochastic dominance}.
In one way or another, they prevent, to some extent, negative outcomes related to the so-called black swan scenarios (i.e., those with high cost and low probability).

In this work we introduce the {\it outsourcing-risk aversion} (ORA) formulation for two-stage DRO as a particular case of the second-order stochastic dominance measure.
It aims to prevent outsourcing-related solutions that incur high costs.
Note that this cost can be balanced, in model DRO-RN \eqref{robust_RN}, with lower costs in other scenarios for the probability distribution $\mathbb{P}$ selected by the model from the ambiguity set $\mathcal{P}$. On the other hand, it can be shown that it is a coherent risk measure, as defined in \cite{Artzner99}, and it is partly inspired by the penalization of the infeasibility of the constraints with stochastic elements that is considered in \cite{KleinHaneveld86}.
 Its peculiarity consists of considering a high cost (usually, due to outsourcing) prevention policy $b$ in the cost profile set $\BB$.
Notice that the ORA formulation can consider multiple high costs $u$ to be minimized jointly with the scenario expected cost surplus penalization of multiple cost $\iota$-thresholds in the ambiguity set. Its advantage is precisely that it allows each policy $b\in \BB$ to protect the DRO solution against modeler-driven outsourcing cost levels.
 To do that, for each policy $b$, the ORA formulation \eqref{robust_RA} requires the incorporation of {\it a)} the variable $u^b$ that is the highest scenario expected cost in the ambiguity set $\P$,  where the binary variable $\gamma_{\mathbb{P}}^b$ has the value 1 \textit{iff} member $\mathbb{P}\in\P$ is that cost's provider, {\it b)} a modeler-driven cost threshold $\iota^b$, {\it c)} the variable $s^{\omega,b}$ as the surplus over cost threshold $\iota^b$ under scenario $\omega$, and {\it d)} the penalization unit $\rho^b$ for the scenario expected cost surplus in ambiguity set member $\mathbb{P}\in\P: \, \gamma_{\mathbb{P}}^b=1$.

 Then, the so-called DRO-ORA metamodel can be expressed as
\begin{subequations}\label{robust_RA}
 \begin{eqnarray}
    && \label{DRO_RA-f} \min
      \sum_{b\in\BB} \kappa^b \bigl( u^b + \rho\sum_{\mathbb{P}\in\P}\sum_{\omega\in\Omega_{\mathbb{P}}} w^{\omega}s^{\omega,b}\gamma_{\mathbb{P}}^b \bigr)\\
   && \mbox{s.to}  \text{ cons system } \eqref{DRO_RN-c}-\eqref{DRO_RN-d}\\
   && \label{DRO_RA-b}
    \underline{u^b}(1-\gamma_{\mathbb{P}}^b) + u^b \gamma_{\mathbb{P}}^b \le
    c^{T}\delta + \sum_{\omega\in\Omega_{\mathbb{P}}}w^{\omega}\nu_1^{\omega T} G^{\omega} \nu_2^{\omega}
    \le u^b                               ~~\forall \mathbb{P}\in \P, \, b\in\BB\\
   && \label{robust_RA-e}
       \bigl ( c^{T}\delta +\sum_{\omega\in\Omega_{\mathbb{P}}}
                        w^{\omega}\nu_1^{\omega T} G^{\omega} \nu_2^{\omega}  - s^{\omega,b} \bigr ) \gamma_{\mathbb{P}}^b \le \iota^b
                                          ~~ \forall \omega\in\Omega_{\mathbb{P}}, ~~\mathbb{P}\in\P, \, b\in\BB \\
   &&\label{robust_RA-g} \sum_{\mathbb{P}\in \P} \gamma_{\mathbb{P}}^b=1 ~~ \forall b\in\BB\\
   &&\label{robust_RA-h} \gamma_{\mathbb{P}}^b\in \{0,1\},  ~ 0\le s^{\omega,b}\gamma_{\mathbb{P}}^b \le \overline{s}^b
                                          ~~ \forall \omega\in\Omega_{\mathbb{P}}, ~~\mathbb{P}\in\P, ~b\in\BB.
    \end{eqnarray}
  \end{subequations}

The objective function \eqref{DRO_RA-f} consists of the weighted (through parameter $\kappa$) highest expected cost $u^b$ co-minimized jointly with the penalization of the expected cost surplus in the scenarios for the ambiguity set member that is the provider of cost $u^b$
(i.e., $\mathbb{P}\in\P: \gamma_{\mathbb{P}}^b=1$), for each policy $b\in\BB$.
Observe that the constraint system \eqref{DRO_RA-b}
is satisfied with strict equality for that ambiguity set $\mathbb{P}$.
Otherwise, the input bound $\underline{u}^b$ and variable bound $u^b$ are satisfied.
Constraint system \eqref{robust_RA-e}-\eqref{robust_RA-h} defines the bounded surplus over cost threshold $\iota^b$ after each realization of $\tilde{\xi}$.
Notice that the {\it special ordered set S1}-based constraint \eqref{robust_RA-g} defines the binary variable $\gamma_{\mathbb{P}}^b=1$ to identify the ambiguity set member that is the highest cost provider for policy $b$.

\subsubsection{Mathematically equivalent linearized reformulation}\label{qp-to-lp}
Given the high difficulty on solving the risk-neutral \eqref{robust_RN} and the risk-averse \eqref{robust_RA} versions of this NP-hard combinatorial problem, a mathematically equivalent reformulation is proposed; it is so-called Linearized mixed Integer Programming (LIP).
So, the replacement of the quadratic expressions is as follows:
\begin{enumerate}
\item The expression $\nu_1^{\omega T} G^{\omega} \nu_2^{\omega}$ in \eqref{DRO_RN-b} and \eqref{robust_RA-e} is replaced with $G'^{\omega}v^\omega$,
    where $G'^{\omega}$ is a row vector and $v^\omega$ is a variables vector for $v^\omega\in[0,1]^{m\times n}$ as in \cite{Guignard12};
    it considers the so-called Reformulation Linearization Technique (RLT1) introduced in \cite{AdamsSherali86},
    where the set of RLT1 inequalities is appended.
    Vector $G'^{\omega}$ is composed of the ordered row vectors  $G_i^\omega\in\mathbb{R}^n \, \forall i=1,..., m$,
    in parameter matrix $G^\omega\in\mathbb{R}^{m,n}$,
    such that ${G'}^\omega=(G_1^\omega, G_2^\omega,..., G_m^\omega)$.
    Note: The optimal RLT1 solution has a binary value for $v^\omega$.

\item The bilinear expressions $u^b \gamma_{\mathbb{P}}^b$ in \eqref{DRO_RA-b} and $s^{\omega,b} \gamma_{\mathbb{P}}^b$ in \eqref{DRO_RA-f} and \eqref{robust_RA-h} are replaced with the continuous variables, say,
    $u_{\mathbb{P}}^b\in\mathbb{R}$ and $s^{\omega,b}_{\mathbb{P}}\in\mathbb{R}$, respectively, and the appending of the related so-called Fortet inequalites,
    see \cite{Fortet60,McCormick76}.

\item The expression
      $\bigl [ c^{T}\delta +\sum_{\omega\in\Omega_{\mathbb{P}}}w^{\omega}\nu_1^{\omega T} G^{\omega} \nu_2^{\omega} - s^{\omega'b} \bigr ] \gamma_{\mathbb{P}}^b$ in \eqref{robust_RA-e}, is first replaced with
       $F^{\omega}\gamma_{\mathbb{P}}^b - s^{\omega}\gamma_{\mathbb{P}}^b$, where
      $F^\omega= c^{T}\delta +\sum_{\omega\in\Omega_{\mathbb{P}}}w^{\omega}\nu_1^{\omega T} G^{\omega} \nu_2^{\omega}$ and, second, $F^{\omega}\gamma_{\mathbb{P}}^b$ itself is replaced with variable, say, $F^{\omega,b}_{\mathbb{P}}\in\mathbb{R}$ and the related Fortet inequalites.
\end{enumerate}
\section{The  CDDP-DRO risk-neutral model and its Scenario Cluster Decomposition}\label{sec:RN}
This section presents three fundamental aspects concerning the
modeling and treatment of the CDDP robust optimization problem.
Subsection \ref{sec:LIP-RN} presents the specific  formulation of the
so-called CDDP-DRO model in its risk-neutral version.
Subsection \ref{sec:LIP-RN-SCD} deals with the {\it scenario cluster decomposition} of the model.
Subsection \ref{sec:rn-lb-ub} introduces the \textit{min-max} matheuristic
scheme for obtaining (hopefully, tight) lower and upper bounds.

\subsection{Linearized risk-neutral CDDP-DRO model (LIP-RN)}\label{sec:LIP-RN}
Based on the previous formulation \eqref{robust_RN} and incorporating the problem-specific notational elements, the proposed two-stage robust risk-neutral framework for the design and operation of a Cross-Dock Door platform so-called CDDP-DRO model can be constructed as follows.

The objective of the problem will include two types of costs to be
minimized: the cost associated with the door's building and installing
their nominal capacities ($C_{1}$), and the operating cost of the
platform ($F^{\omega}$), which will take a different value under each scenario for the ambiguity set.

The CDDP-DRO risk-neutral model can be expressed as
\begin{subequations}\label{LIP-RN}
\bea
&& \label{LIP-RN-1}  z_{RN}^* \, = \, min \, u\\
&& \label{LIP-RN-2} \displaystyle \text{s.to }  C_1 + \sum_{\omega\in\Omega_\mathbb{P}} w^\omega F^\omega \leq u ~~\forall \mathbb{P}\in\P\\
&& \label{LIP-RN-3}  \text{1st and 2nd stage cons systems } \eqref{dro-1st-c} \text{ and } \eqref{dro-2nd-c}, \text{ resp.}
\eea
\end{subequations}

\subsubsection{First-stage constraints}
The first-stage constraint system \eqref{dro-1st-c} can be expressed as
\begin{subequations}\label{dro-1st-c}
\bea
&& \label{Q1}  \displaystyle
\alpha_{ki}\in\{0,1\} \,\, \forall k\in\K_i, \,\, \sum_{k\in\K_i}\alpha_{ki} \leq 1 \,\, \forall i\in\I, \,\,
\sum_{i\in\I}\sum_{k\in\K_i}\alpha_{ki} \leq \overline{I}\\
&& \label{Q2}  \displaystyle
\beta_{kj}\in\{0,1\} \,\,  \forall k\in\K_j, \,\, \sum_{k\in\K_j}\beta_{kj} \leq 1   \,\, \forall j\in\J, \,\,
\sum_{j\in\J}\sum_{k\in\K_j} \beta_{kj} \leq \overline{J}\\
&& \label{C1} \displaystyle
 C_1 \, = \, \sum_{i\in\I}\sum_{k\in\K_i}F_{ki} \alpha_{ki} + \sum_{j\in\J}\sum_{k\in\K_j}F_{kj} \beta_{kj}.
 \eea
\end{subequations}
 where $\alpha_{ki} \, (\beta_{kj})$ denote the first stage variables
that take the value 1 iff capacity level $k$ is installed in strip door $i$ (stack door $j$),
 for $k\in\K_i$ ($k\in\K_j$) and, 0 otherwise. $\I$ ($\J$) denotes the strip  (stack) door candidate set,
 without including inbound  outsourcing `door' $i=0$ ($j=0$),
  and $\overline{I}$ ($\overline{J}$) denotes the  upper bound on the number of strip (stack) doors
 in the CD platform infrastructure.
System \eqref{dro-1st-c} defines the doors' building, nominal capacity installation and
related cost $C_1$, where $F_{ki}$ ($F_{kj}$) is the  cost of
installing capacity level $k$ in strip door $i$ (stack door $j$), for
$k\in\K_i$ ($k\in\K_j$), and $i\in\I$ ($j\in\J$).
\subsubsection{Second-stage constraints}
The second stage constraint system \eqref{dro-2nd-c} defines the node--door
standard and outsourcing assignment constraining in the CD platform
and the related cost $F^{\omega}$ under the scenarios for each ambiguity set member.
It can be expressed as
\begin{subequations}\label{dro-2nd-c}
\bea
&& \nonumber \text{For } \omega\in\Omega_{\mathbb{P}}, ~\mathbb{P}\in\P:\\
&& \label{Q3} \displaystyle
x_{mi}^\omega\in\{0,1\} \,\, \forall i\in\I_m^\omega\cup\{0\}, \,\,
\alpha_0^\omega\in\{0,1\}, \,\,
x_{m0}^\omega \leq \alpha_0^\omega, \,\
\sum_{i\in\I_m^\omega\cup\{0\}} x_{mi}^\omega = 1 \,\, \forall m\in\M^\omega\\
&& \label{Q4} \displaystyle
\sum_{m\in\M^\omega:i\in\I_m^\omega}S_m^\omega x_{mi}^\omega \leq (1- D_i^\omega)\sum_{k\in\K_i} S_{ki} \alpha_{ki} \,\, \forall i\in\I\\
&& \label{Q5} \displaystyle
y_{nj}^\omega\in\{0,1\} \,\, \forall j\in\J_n^\omega\cup\{0\}, \,\,
\beta_0^\omega\in\{0,1\}, \,\,
y_{n0}^\omega \leq \beta_0^\omega, \,\,
\sum_{j\in\J_n^\omega\cup\{0\}} y_{nj}^\omega = 1 \,\, \forall n\in\N^\omega\\
&& \label{Q6} \displaystyle
\sum_{n\in\N^\omega:j\in\J_n^\omega}R_n^\omega y_{nj}^\omega \leq (1- D_j^\omega)\sum_{k\in\K_j} R_{kj}\beta_{kj} \,\, \forall j\in\J\\                                                                     && \label{F-omega} \displaystyle
F^\omega \, = \,  F_0(\alpha_0^\omega + \beta_0^\omega)
                + \sum_{m\in\M^\omega}\sum_{i\in\I_m^\omega\cup\{0\}}\sum_{n\in\N^\omega}\sum_{j\in\J_n^\omega\cup\{0\}}
                  G_{minj}^\omega v_{minj}^\omega,
\eea
\end{subequations}
where  the variable $x_{mi}^\omega   = 1$ ($y_{nj}^\omega   = 1$ )  iff
origin node $m$ (destination node $n$) is assigned to strip door $i$
(stack door $j$), and 0 otherwise. Meanwhile, the variable
$x_{m0}^\omega  \, (y_{m0}^\omega ) = 1$  iff node $m$ is \textit{not}
assigned to any strip (stack) door, and 0 otherwise. Likewise, the
variable $\alpha_0^\omega \, (\beta_0^\omega) = 1$ iff inbound
(outbound) "outsourcing door" $i=0$ ($j=0$) is considered, and 0 otherwise.
 Constraints  \eqref{Q3} and \eqref{Q5} force the inbound and outbound outsourcing to be activated
if the assignment of an origin node to any strip door and a destination node to any stack doors is not feasible for a given scenario, respectively.
Constraints \eqref{Q4} and \eqref{Q6} upper bound the commodity volume
handling in the strip and stack doors, based on the doors' net
capacity, where $S_m^\omega$ ($R_n^\omega$) denotes the  inbound (outbound)
volume  entering (exiting) through any strip (stack) door from node
$m$ (at node $n$), and $\M^\omega$
($\N^\omega$) is the origin (destination) node set under scenario
$\omega$. $D_i^\omega$ ($D_j^\omega$) is the disruption fraction of
the nominal capacities, say $S_{ki}$ ($R_{kj}$), of strip door $i$
(stack one $j$).
Constraint \eqref{F-omega} defines the operating cost of the
CD platform ($F^{\omega}$) in terms of the  standard handling and
transportation cost of commodity volume ($G_{minj}^{\omega}$) due to entering the platform  from node $m$ through strip door $i$ and
 exiting to node $n$ through stack door $j$, as well as the  high
 enough penalizations ($F_{0}$, $G_{m0n0}^\omega$) for using the `outsourcing doors' $i=0$ and $j=0$.
The continuous variable $v_{minj}^\omega$ in \eqref{F-omega} (that takes the values either 0 or 1 in the optimal solution) has replaced
the binary  bilinear expression $x_{mi}^\omega y_{nj}^\omega$ as presented in Subsection \ref{qp-to-lp} (bullet 1),
in a mathematically equivalent form, by additionally appending the RLT1 inequalities to be expressed as
\begin{equation}\label{v-lip}
\begin{array}{lll}
&& \text{For } \omega\in\Omega_{\mathbb{P}}, ~~\mathbb{P}\in\P:\\
&& \displaystyle \sum_{j\in\J_n^\omega\cup\{0\}} v_{minj}^\omega = x_{mi}^\omega
      \quad \forall i\in\I_m^\omega\cup\{0\}, \, m\in\M^\omega, \, n\in\N^\omega\\
&& \displaystyle \sum_{i\in\I_m^\omega\cup\{0\}} v_{minj}^\omega = y_{nj}^\omega
      \quad \forall m\in\M^\omega, \, j\in\J_n^\omega\cup\{0\}, \, n\in\N^\omega\\
&& v_{minj}^\omega \geq 0
      \quad \forall i\in\I_m^\omega\cup\{0\}, \, m\in\M^\omega, \, j\in\J_n^\omega\cup\{0\}, \, n\in\N^\omega.
\end{array}
\end{equation}

\subsection{Scenario-cluster-based model decomposition}\label{sec:LIP-RN-SCD}
Our prior experience in solving the primary deterministic model (CDAP), see \cite{cdap24}, together with the results presented in Section \ref{sec:results}, indicates that it is not possible to achieve optimality for model \eqref{LIP-RN} within a reasonable time, even for medium-sized instances.
 Moreover, notice that the members of the ambiguity set $\P$ are only linked with variable $u$ in the model, besides the first stage constraints \eqref{dro-1st-c};
the overall scheme for problem-solving proposed in this work takes advantage of this special structure,
seemingly amenable for  a scenario cluster decomposition.

Therefore, a split-variable reformulation of  model \eqref{LIP-RN} is proposed that enables its decomposition into  smaller submodels related to  clusters of scenarios for each member of the ambiguity set.  This proposal  extends the scenario clustering scheme presented in \cite{cddp-ts24}.

Following \cite{mcld16}, for each cluster $c\in\C$, the {\it splitting variable constraint
system} can be expressed as
\begin{subequations}\label{svc-rn}
\bea
& \label{svc-rn-alpha}  \alpha_{ki}^c -\alpha_{ki}^{n(c)} \leq 0  & \forall k\in\K_i, \, i\in\I\\
& \label{svc-rn-beta}   \beta_{kj}^c -\beta_{kj}^{n(c)} \leq 0    & \forall k\in\K_j, \, j\in\J\\
& \label{svc-rn-u}      u^c -u^{n(c)} \leq 0.
\eea
\end{subequations}
where $\C$  is the total (lexicographically ordered) set of  scenario
clusters, to be expressed as
$\C=\bigcup_{\mathbb{P}\in\P}\C_{\mathbb{P}}$, $\C_{\mathbb{P}}$  is the set of scenario clusters in ambiguity set
member $\mathbb{P}$, and  $\Omega^c$ denotes the scenario cluster $c$,
for $c\in\C_{\mathbb{P}}$.

Note: $\Omega_{\mathbb{P}} = \bigcup_{c\in\C_{\mathbb{P}}}\Omega^c$, where
$\Omega_{\mathbb{P}}$ is the lexicographically ordered set of
scenarios for member $\mathbb{P}$,
and $n(c)$ denotes the next one to cluster $c$ in set $\C_{\mathbb{P}}$.

Under this framework, given a cluster partition of $\Omega_{\mathbb{P}}$ and introducing a
replica  for each cluster $c$ of the variables
$\alpha_{ki}^c$, $\beta_{kj}^{c}$, and $u^{c}$, we derive a mathematically
equivalent reformulation of model \eqref{LIP-RN}, which can be expressed as
\begin{subequations}\label{dro-svc-rn}
\bea
&& \label{dro-of-rn} \displaystyle
z_{RN}^* \, = \, min \, u\\
&& \label{dro-u-rn} \displaystyle
   \text{s.to } C_1
             + \sum_{c\in\C_{\mathbb{P}}} \tilde{w}^c \sum_{\omega\in\Omega^c} \tilde{w}^\omega F^\omega \leq u ~~ \forall \mathbb{P}\in\P\\
&& \text{splitting variable cons
system } \eqref{svc-rn}\\
&& \nonumber \text{cons system } \eqref{dro-1st-c}\text{-}\eqref{dro-2nd-c},
\text{where } \alpha_{ki}, \, \beta_{kj},  \,u, \, C_1 \text{ replaced with } \alpha_{ki}^c, \, \beta_{kj}^c,  \,u^c, \, C_1^c ~~ \forall c\in\C_{\mathbb{P}}, \mathbb{P}\in\P,
\eea
\end{subequations}
where $\tilde{w}^c$ is the weight associated with cluster $c$,
expressed as $\sum_{\omega\in\Omega^c}w^\omega$, for $c\in\C_{\mathbb{P}}$, and,
then                   $\sum_{c\in\C_{\mathbb{P}}}\tilde{w}^c=1$, and
$\tilde{w}^\omega$  denotes the  weight associated with scenario
$\omega$ that belongs to cluster $c$, expressed as
$\frac{w^\omega}{\tilde{w}^c}$ and, then,
$\sum_{\omega\in\Omega^{c}}\tilde{w}^\omega=1$.

It is worth noting that, given \eqref{svc-rn},
it results that the variables $(.)^c$ and $(.)^{c'}$ have the same
value for clusters $c,c'\in\C_\mathbb{P}$.
Similarly, $u^{\overline{c}}$ has the same value as $u$ \eqref{LIP-RN-2}.
Anyway, model \eqref{dro-svc-rn} makes it possible to design a decomposition approach, as described below.

The relaxation of the  constraint system \eqref{svc-rn} results in
independent $c$-submodels \eqref{RN-SCD-c},  for each cluster $c\in\C_\mathbb{P}$, to be expressed as
\begin{subequations}\label{RN-SCD-c}
\bea
&& z_{RN}^c \, = \, min \, u^c \\
&& \displaystyle
   \text{s.to } C_1^c + \sum_{\omega\in\Omega^c} \tilde{w}^\omega F^\omega \leq u^c \\
&& \text{cons system } \eqref{dro-1st-c}\text{-}\eqref{dro-2nd-c},
   \text{where } \alpha_{ki}, \, \beta_{kj}, \,u, \, C_1 \text{ replaced with } \alpha_{ki}^c, \, \beta_{kj}, \,u^c, \, C_1^c.
\eea
\end{subequations}

\subsection{Obtaining lower ($\underline{z}_{RN}$) and upper ($\overline{z}_{RN}$) bounds  in model  \eqref{LIP-RN}} \label{sec:rn-lb-ub}
There are a myriad of decomposition algorithms for solving general two-stage stochastic MILP models,
see, e.g., \cite{mcld16} for an overview of the literature. However, given the inherent computational difficulty of solving even the scenario submodels in this problem, any iterative strategy must be ruled out.
Taking this limitation into account, and exploiting the cluster-based decomposition of the model,  we propose below a two-step \textit{min-max} {\it matheuristic} procedure that yields lower and upper bounds for the problem with a reasonable computational effort. On the one hand, this approach enables the derivation of a lower bound for the objective function, together with a set of candidate solutions for the first stage. On the other hand, once these candidate solutions have been evaluated, a suitable feasible solution is obtained, thereby yielding an upper bound for the problem.

The main steps of the scheme are as follows:
\begin{enumerate}\parskip 0.5mm
\item Obtaining the lower bound $\underline{z}_{RN}$:
\begin{enumerate}\parskip 0.5mm
\item Solve the $c$-submodel \eqref{RN-SCD-c} and retrieve the value of the related
  objective function $z_{RN-SCD}^c$  and vector $( \hat{\alpha}^c, \, \hat{\beta}^c )$ of the first-stage variables $(\alpha_{ki}^c \, \forall ki, \, \beta_{kj}^c \, \forall kj)$  for each cluster $c\in\C$.

\item Solve the LP submodel in the B\&C root node for solving model \eqref{LIP-RN}, 
 and retrieve the value of the related objective function, say, $z_{LP-RN}$.

\item The \textit{max-max} (lower) bound $\underline{z}_{RN}$ can be expressed as
\begin{equation}\label{lb-RN}\displaystyle
\underline{z}_{RN} \, = \, max \{ z_{LP-RN}, \, max_{\mathbb{P}\in\P}\{ \sum_{c\in\C_{\mathbb{P}}} \tilde{w}^{c} z_{RN}^c \} \}.
\end{equation}
\end{enumerate}

\item Obtaining the upper bound $\overline{z}_{RN}$:
\begin{enumerate}\parskip 0.5mm
\item Let $( \hat{\alpha}^c, \, \hat{\beta}^c )$ be the vector of the first-stage variables $(\alpha_{ki}^c \, \forall ki, \, \beta_{kj}^c \, \forall kj)$ from the solution of the $c$-submodel \eqref{RN-SCD-c} as obtained in Step 1(a), and
      compute the first stage cost, say, $\hat{C}_1^c$ \eqref{C1},  for each cluster $c\in\C$.

\item For each $\omega\in\Omega_{\mathbb{P}}, ~\mathbb{P}\in\P$, solve the $\omega$-submodels (\ref{RN-cdap}), see \cite{cdap24}, where the vector of the $(\alpha,\beta)$-first stage variables takes the solution
      $( \hat{\alpha}^c, \, \hat{\beta}^c )$, for  each cluster  $c\in\C_\mathbb{P}$.
\begin{equation}\label{RN-cdap}
\begin{array}{lll}
&& z_{RN}^{c,\omega} \, = \, min \, F^{c,\omega} \eqref{F-omega}\\
&& \text{s. to cons system } \eqref{dro-1st-c}\text{-}\eqref{dro-2nd-c}, \text{ where only scenario } \omega \text{ is considered}
 \end{array}
\end{equation}

\item The \textit{min-max} (upper) bound $\overline{z}_{RN}$ can be expressed as
\begin{equation} \label{ub-RN} \displaystyle
\overline{z}_{RN} \, = \,  min_{c\in\C}\{ \hat{C}_1^c + max_{\mathbb{P}\in\P} \{ \sum_{\omega\in\Omega_\mathbb{P}} w^\omega z_{RN}^{c,\omega} \} \}.
\end{equation}
\end{enumerate}
\end{enumerate}
\section{The  CDDP-DRO risk-averse  model and its Scenario Cluster Decomposition}\label{sec:ORA}
As a counterpart of the DRO risk-neutral model \eqref{LIP-RN}, the CDDP-DRO risk-averse version and its related scenario cluster decomposition are now introduced as well as the
\textit{min-max} matheuristic for obtaining lower and upper bounds for the solution value.
It can be observed that the ORA metamodel \eqref{robust_RA}  preserves in essence the original structure of the risk-neutral formulation \eqref{robust_RN} in order to perform its scenario cluster decomposition. Since a singleton policy profile (i.e., $|\BB|=1$) is sufficient to illustrate this decomposition, it can be adopted without loss of generality for exposition purposes. So, throughout the remainder of the paper, the superscript $b$ is omitted from the corresponding parameters and variables.

The section is organized as done for Section \ref{sec:RN}.
\subsection{Linearized Risk-Averse CDDP-DRO Model (LIP-ORA)}\label{sec:LIP-ORA}
The system for determining $u$ in the linearized reformulation of the  ORA metamodel can be expressed as
\begin{subequations}\label{DRO-ORA-u}
\bea
&& \displaystyle
   \underline{u} (1-\gamma_{\mathbb{P}}) + u_{\mathbb{P}} \leq \sum_{\omega\in\Omega_{\mathbb{P}}}w^\omega C_{12,\mathbb{P}}^\omega \leq u
                                                                 ~~ \forall \mathbb{P}\in\P\\
&& \label{dro-gamma} \displaystyle
   \sum_{\mathbb{P}\in\P} \gamma_{\mathbb{P}} = 1 \\
&& \label{dro-up-Fortet}
   u_{\mathbb{P}} \leq \overline{u} \gamma_{\mathbb{P}}, \, u_{\mathbb{P}} \leq u, \,
   u \leq u_{\mathbb{P}} + \overline{u} (1- \gamma_{\mathbb{P}}) ~~ \forall \gamma_{\mathbb{P}}\in\{0,1\}, ~\mathbb{P}\in\P,
\eea
\end{subequations}
where $u_{\mathbb{P}}$ is a continuous variable that replaces the bilinear term  $u \gamma_{\mathbb{P}}$,
$\underline{u}$ and $\overline{u}$ are the highest lower bound and the smallest upper bound of variable $u$ that are allowed, respectively,
and the Fortet inequalities \eqref{dro-up-Fortet} follow Subsection \ref{qp-to-lp} (bullet 2).
Also, inspired in metamodel \eqref{robust_RA}, the  constraint system for preventing outsourcing-based solutions can be expressed as
\begin{subequations}\label{DRO-ORA}
\bea
&& \nonumber \text{For } \omega\in\Omega_{\mathbb{P}}, ~~ \mathbb{P}\in\P:\\
&& \label{DRO-ORA-0}
  C_{12}^\omega = C_1 \eqref{C1} + F^\omega \eqref{F-omega}\\
&& \label{DRO-ORA-1}  \displaystyle
  C_{12,{\mathbb{P}}}^\omega - s_{\mathbb{P}}^\omega \leq \iota, \quad  0 \leq s_{\mathbb{P}}^\omega \leq \overline{s}\\
&& \label{DRO-ORA-Fortet-C}
C_{12,{\mathbb{P}}}^\omega \leq \overline{C}_{12} \gamma_{\mathbb{P}}, \,
C_{12,{\mathbb{P}}}^\omega \leq C_{12}^\omega, \, C_{12}^\omega \leq C_{12,{\mathbb{P}}}^\omega + \overline{C}_{12} (1- \gamma_{\mathbb{P}})\\
&& \label{DRO-ORA-Fortet-s}
s_{\mathbb{P}}^\omega \leq \overline{s} \gamma_{\mathbb{P}}, \,
s_{\mathbb{P}}^\omega \leq s^\omega, \, s^\omega \leq s_{\mathbb{P}}^\omega + \overline{s} (1- \gamma_{\mathbb{P}}),
\eea
\end{subequations}
where $C_{12}^\omega$ is the cost composed of the first stage one and second stage cost under scenario $\omega$, in ambiguity set member ${\mathbb{P}}$, so that  $C_{12}^\omega$ = $C_1$ + $F^\omega$;
$\overline{C}_{12}$ is the modeler-driven upper bound on $C_{12}^\omega$;
$C_{12,{\mathbb{P}}}^\omega$ is the continuous variable given the cost under scenario $\omega$ in the definition of $u$ \eqref{DRO-ORA-u},
\textit{iff} it is provided by member $\mathbb{P}\in\P:\gamma_{\mathbb{P}}=1$; otherwise, 0, so that
$C_{12,{\mathbb{P}}}^\omega = C_{12}^\omega \gamma_{\mathbb{P}}$; and
$s_{\mathbb{P}}^\omega$ is the surplus variable over cost threshold $\iota$ under scenario $\omega$, for $\gamma_{\mathbb{P}}=1$, so that
$s_{\mathbb{P}}^\omega = s^\omega \gamma_{\mathbb{P}}$. Constraints \eqref{DRO-ORA-1} define the bounded cost surplus over the threshold $\iota$ in any scenario.
And the Fortet inequalities \eqref{DRO-ORA-Fortet-C} and \eqref{DRO-ORA-Fortet-s} follow Subsection \ref{qp-to-lp} (bullets 2 and 3).
 So, model LIP-ORA can be expressed as
\begin{equation}\label{LIP-ORA}
\begin{array}{lll}
&&  \displaystyle
z_{ORA}^*\, = \, min \, u  + \rho \sum_{\mathbb{P}\in\P}\sum_{\omega\in\Omega_{\mathbb{P}}} w^{\omega}s_{\mathbb{P}}^{\omega}\\
&& \text{ s.to cons system } \eqref{dro-1st-c},\eqref{dro-2nd-c},\eqref{DRO-ORA-u},\eqref{DRO-ORA}.
\end{array}
\end{equation}
Notice that $u_{\mathbb{P}}=C_{12,{\mathbb{P}}}^\omega=0$  for $\omega\in\Omega_{\mathbb{P}}, \, \mathbb{P}\in\P: \gamma_{\mathbb{P}}=0$, see S1 \eqref{dro-gamma}.
\subsection{Scenario cluster decomposition of model LIP-ORA}\label{sec:LIP-ORA-SCD}
A high computational effort may be required for solving model LIP-ORA \eqref{LIP-ORA} even for medium-sized instances,
being, in general, higher than the effort required for solving model LIP-RN \eqref{LIP-RN}.
Moreover, notice that, besides the first stage constraints \eqref{Q1}-\eqref{Q2},
the members of the ambiguity set $\P$ are linked in the model with constraint \eqref{dro-gamma}.

The overall scheme proposed in this work exploits the model’s special structure, which is amenable to decomposition.
 To this end, let $\gamma_{\mathbb{P}}, \, u_{\mathbb{P}}, \, u$ from \eqref{DRO-ORA-u} be replaced with $\gamma^c, \, u_c, \, u^c, \,\, c\in\C_{\mathbb{P}}, \, \mathbb{P}\in\P$.
So, the new splitting variable constraint system can be expressed as
\begin{equation}\label{svc-ORA}
\begin{array}{lll}
&& \text{cons system } \eqref{svc-rn}  ~~\forall c\in\C_{\mathbb{P}}\\
&& \gamma^c -\gamma^{n(c)} \leq 0      ~~\forall c\in\C_{\mathbb{P}}, ~~ \mathbb{P}\in\P.
\end{array}
\end{equation}

Based on \eqref{svc-ORA}, an equivalent formulation to LIP-ORA \eqref{LIP-ORA}, can be expressed as
\begin{equation}\label{dro-svc-ORA-cons}
\begin{array}{lll}
&& \displaystyle
z_{ORA}^* \, = \, min \, u
                      + \rho \sum_{\mathbb{P}\in\P}\sum_{c\in\C_{\mathbb{P}}}\tilde{w}^c\sum_{\omega\in\Omega^c} \tilde{w}^\omega s_{\mathbb{P}}^\omega\\
&& \displaystyle
 \text{s.to }
 \underline{u} (1-\gamma_{\mathbb{P}}) + u_c \leq
 \sum_{c\in\C_{\mathbb{P}}}\tilde{w}^c\sum_{\omega\in\Omega^c}\tilde{w}^\omega C_{12,\mathbb{P}}^\omega
                                                     \leq u   ~~\forall \mathbb{P}\in\P\\
&&  u_c \leq \overline{u} \gamma^c, \, u_c \leq u^c , \,
   u^c  \leq u_c + \overline{u} (1- \gamma^c), \quad \gamma^c\in\{0,1\}   ~~\forall c\in\C_{\mathbb{P}}, \, \mathbb{P}\in\P\\
&& \text{cons system } \eqref{dro-1st-c}\text{-}\eqref{dro-2nd-c},\eqref{dro-gamma},\eqref{DRO-ORA},\eqref{svc-ORA}, \text{ where }\\
&& \gamma_{\mathbb{P}}, \, \Omega_{\mathbb{P}}  \text{ replaced with }
   \gamma^c, \, \Omega^c,    ~~\forall c\in\C_{\mathbb{P}}, ~~ \mathbb{P}\in\P, \\
&& \alpha_{ki}, \, \beta_{kj} , \,u, \, C_1 \text{ replaced with }
             \alpha_{ki}^c, \, \beta_{kj}^c, \,u^c , \, C_1^c \text{ for } c\in\C_{\mathbb{P}}, ~~ \mathbb{P}\in\P.
\end{array}
\end{equation}

The relaxation of the constraint system \eqref{dro-gamma}, \eqref{svc-ORA} results in  the  $c$-submodels \eqref{ORA-SCD-c} that can be expressed as
\bcb
\begin{equation}\label{ORA-SCD-c}
\begin{array}{lll}
&& \text{For cluster } c\in\C_{\mathbb{P}}:\\
&& \displaystyle
z_{ORA}^c \, = \, min \, u^{c} + \rho \sum_{\omega\in\Omega^c} \tilde{w}^\omega s_{\mathbb{P}}^\omega\\
&& \displaystyle
 \text{s.to } \underline{u} \leq \sum_{\omega\in\Omega^c}\tilde{w}^\omega  C_{12,\mathbb{P}}^\omega \leq u^c\\
&& \text{cons system } \eqref{dro-1st-c}\text{-}\eqref{dro-2nd-c}, \eqref{DRO-ORA},
   \text{ where } \gamma_{\mathbb{P}}=\gamma^c=1 \text{ is assumed, and }\\
&& \alpha_{ki}, \, \beta_{kj}, \, C_1, \, \Omega_{\mathbb{P}},  \text{ replaced with }
   \alpha_{ki}^{c}, \, \beta_{kj}^{c}, \, C_1^{c}, \, \Omega^{c}.
\end{array}
\end{equation}
\subsection{Obtaining the lower ($\underline{z}_{ORA}$) and upper ($\overline{z}_{ORA}$) bounds  of $z_{ORA}^*$ in model LIP-ORA \eqref{LIP-ORA}}\label{sec:ORA-lb-ub}
The scheme follows the pattern considered for obtaining the incumbent solution in the  risk-neutral version, see Section \ref{sec:rn-lb-ub}.
\begin{enumerate}\parskip 0.5mm
\item Obtaining the lower bound $\underline{z}_{ORA}$:
\begin{enumerate}\parskip 0.5mm
\item Solve the cluster $c$-submodel \eqref{ORA-SCD-c} and retrieve its
  solution value $z_{ORA}^c$  and vector $( \hat{\alpha}^c, \, \hat{\beta}^c )$ of the first-stage variables, for $c\in\C$.

\item Solve the LP relaxation in the B\&C root node for solving model \eqref{LIP-ORA}, 
 and retrieve its solution value, say, $z_{LP-ORA}$.

\item The \textit{max-max} bound $\underline{z}_{ORA}$ can be expressed as
\begin{equation}\label{lb-ORA}\displaystyle
\underline{z}_{ORA} \, = \, max \{ z_{LP-ORA}, \, max_{\mathbb{P}\in\P}\{ \sum_{c\in\C_{\mathbb{P}}} \tilde{w}^{c} z_{ORA}^c \} \}.
\end{equation}
\end{enumerate}

\item Obtaining the upper bound $\overline{z}_{ORA}$:
\begin{enumerate}\parskip 0.5mm
\item Let $( \hat{\alpha}^c, \, \hat{\beta}^c )$ as obtained in Step 1(a), and
      compute the first stage cost, say, $\hat{C}_1^c$ \eqref{C1}, for $c\in\C$.

\item For each $\omega\in\Omega_{\mathbb{P}}, ~~ \mathbb{P}\in\P$, solve the scenario $\omega$-submodels `extended'  (\ref{ORA-cdap}), where the vector of the $(\alpha,\beta)$-first stage variables takes the solution
      $( \hat{\alpha}^c, \, \hat{\beta}^c )$,  for $c\in\C_\mathbb{P}$.
\begin{equation}\label{ORA-cdap}
\begin{array}{lll}
&& z_{ORA}^{c,\omega} \, = \, min \, F^{c,\omega} \eqref{F-omega} + \rho s_{\mathbb{P}}^\omega\\
&& \text{s.to cons system } \eqref{dro-1st-c}\text{-}\eqref{dro-2nd-c}, \text{ where only scenario } \omega \text{ is considered } \\
&& F^{c,\omega} - s_{\mathbb{P}}^\omega \leq \iota - \hat{C}_1^c, \,\, 0 \leq s_{\mathbb{P}}^\omega \leq \overline{s}.
\end{array}
\end{equation}

\item The \textit{min-max} bound $\overline{z}_{ORA}$ can be expressed as
\begin{equation} \label{ub-ORA} \displaystyle
\overline{z}_{ORA} \, = \,  min_{c\in\C}\{ \hat{C}_1^c +
max_{\mathbb{P}\in\P} \{ \sum_{\omega\in\Omega_{\mathbb{P}}} w^\omega
{\hat F}^{c,\omega}  \} \}. 
\end{equation}
\end{enumerate}
\end{enumerate}

It is worth remarking that there is a high computational complexity derived from  the special structure of models LIP-RN \eqref{LIP-RN} and LIP-ORA \eqref{LIP-ORA}.
Additionally, the scenario submodels \eqref{RN-cdap} and their ‘extended’ counterparts \eqref{ORA-cdap} exhibit high complexity for each scenario in the ambiguity set once the first-stage variables are fixed.
Consequently, the second stage of the proposed
\textit{min-max} scheme incorporates the recently developed ad hoc constructive
matheuristic introduced in \cite{cdap24}.

As a matter of fact, any iterative scheme for the CDDP problem solving where a set of scenario-related submodels are to be considered,
can not pay the effort.
Observe that those models should have, \textit{at each iteration}, a scenario cluster (or a scenario itself)-related CDDP type of constraint system jointly with a risk-averse one.
Its solving is a very difficult task as it can be deduced from the results of the computational experiment that is reported in this work,
where the proposed non-iterative procedure is considered.
As an example, let the Lagrangean Decomposition (LD) of model CDAP, a non-DRO subset of CDDP that is considered in \cite{cdap24}.
Its problem-solving process is so complex that any LD iteration does require a high computational effort
and, additionally, the weak increase in the LD's lower bound does not pay the effort.
        \section{Computational experiment}\label{sec:results}
The section presents the results of a computational experiment carried out  in order to show the effectiveness of the proposed approach, where the testbed  is composed of versions of the instances I1, I3 and I7, originally introduced in
\cite{cddp-ts24}.
 Additionally, a comparative study was performed between the straightforward use of the state-of-the-art solvers CPLEX and Gurobi.

Appendix B in the Supplementary file introduces a procedure for generating a candidate ambiguity set and deals with obtaining the Wasserstein ambiguity set $\P$.
With regard to the uncertain parameters that appear in the second stage of the CDDP problem, the functional relationships that exist among some of them must be taken into account.
Thus, the parameters that represent the inbound volume entering the platform and the outbound volume exiting from, say $S_m^\omega$ and $R_n^\omega$, respectively, and the commodity handling and transportation cost, say $G_{minj}^\omega$, are directly obtained from $H_{mn}^\omega$, the commodity volume to be consolidated in the platform from origin node $m$ to destination node $n$.
Then, the uncertain (independent) parameters to which the proposed methodology has been applied are: {\it a)} the commodity volume $H_{mn}^\omega$, and    {\it b)} the  capacity disruption fraction of the strip and stack doors $D_i^\omega$ and $D_j^\omega$, respectively.
Taking into account the stochasticity of these second stage parameters, different probability distributions were considered to generate the ambiguity set. On the one hand, the Normal, Weibull, Gamma and Lognormal  distributions are used for generating the different perturbations of the realizations $H_{mn}^\omega$.  On the other hand, the perturbations of the  capacity disruptions $D_i^\omega$ and $D_j^\omega$ were generated using the Uniform distribution.

The experiment has been performed for a singleton policy profile (i.e., $|\BB|=1$) related to the second-order stochastic dominance risk averse measure.
Appendix C in the Supplementary file provides a detailed description of the data in the nominal distribution and the ambiguity sets that have been generated for different versions of each of the three instances considered in the experiment.
The appendix further provides an account of additional details retrieved from the solutions.


\noindent\textbf{Computational settings}. The experiment reported  was conducted on a
Debian GNU/Linux 12 Work Station,  with ES-2670 v2 processors (20 physical cores, 40 logical processors, using up to 20 threads), at 2.50 Ghz and with 128 GB of RAM.
The CPLEX (v22.1.1) Concert Technology library and Gurobi (v13.0.1) were used, embedded in a C/C++  experimental code,
for solving the decomposition submodels as well as the straightforward solving of the full models;
the default parameters were used unless explicity mentioned.

\subsection{Instance I1}\label{sec:I1}
The smallest instance in the experiment is composed of up to 8 origins and 8 destinations, and 4 strip doors and 4 stack doors.
The ambiguity sets $\P$ under consideration are, say, $\P_{5\%}$ and $\P_{10\%}$, where $|\P_{5\%}|=4$ and $|\P_{10\%}|=8$ and $\P_{5\%}\subset\P_{10\%}$;
see Appendix C.1 for the sets generation.

\subsubsection{I1. On solving models LIP-RN \eqref{LIP-RN} and LIP-ORA \eqref{LIP-ORA} for ambiguity set $\P_{5\%}$}\label{sec:I1-RN-P5}
Table \ref{I1.dim-results} shows the dimensions of models LIP-RN \eqref{LIP-RN} and LIP-ORA \eqref{LIP-ORA} as well as the solution cost and the wall time required by the straightforward use of the solvers CPLEX and Gurobi.
The headings are as follows.
$Inst$-$|\P|$, instance code, where the cardinality  of the ambiguity set is included;
$m$, $n01$, $nc$ and $nz$, number of constraints, binary variables, continuous variables and nonzero constraint matrix elements, respectively,
$z^*$, cost of the optimal solution;
$t_{cpx}$ and $t_{grb}$, wall time (seconds) required by the solvers to obtain $z^*$.

\begin{table}[h!]
\begin{center}
{\small
\caption{Instance I1. CDDP models LIP-RN \eqref{LIP-RN} and LIP-ORA \eqref{LIP-ORA}: Dimensions and results of CPLEX and Gurobi \\
Note: $*$: 0.03\% optimality gap}\label{I1.dim-results}
\begin{tabular}{l|cccc|c|cc}
\hline
$Inst$-$|\P|$ & $m$ & $n01$ & $nc$ & $nz$ & $z^*$ & $t_{cpx}$ & $t_{grb}$ \\
\hline
I1-RN-4 &13634 &1680 &32021 &98398  &7260.56 &   856 &   353\\
I1-ORA-4 &13735 &1684 &32065 &99606  &7459.31 &   712 &  3062\\
I1-RN-8 &27258 &3320 &64041 &196717 &7260.56 & 37739 &  6225\\
I1-ORA-8 &27459 &3328 &64129 &199138 &7462.36 &  6365 & 42000$*$\\
\hline
\end{tabular}
}
\end{center}
\end{table}

The optimal solution cost of model LIP-RN \eqref{LIP-RN}, denoted by  $z_{RN}^*= 7260.56$, shown in Table
\ref{I1.dim-results}, is composed of the cost  of building the infrastructure  of the first stage platform
$\hat{C} = 1538$ \eqref{C1} and the second stage expected cost $5722.56$. There are six black swan scenarios in which outsourcing is required.
Table C.3 in the Supplementary file shows the weights and the cross-docking node--door assignment cost $\hat{F}^\omega$ \eqref{F-omega} under the scenarios $\{\omega\}$ in the ambiguity set members.

The aim  of model LIP-ORA \eqref{LIP-ORA} is to prevent the outsourcing costs as identified in the optimal solution of model LIP-RN \eqref{LIP-RN}.
To achieve this, an appropriate modeler-driven triplex (cost threshold, upper bound cost surplus, surplus penalization) is considered in model LIP-ORA \eqref{LIP-ORA}. The optimal solution cost is $z_{ORA}^*= 7459.31$, which is higher than $z_{RN}^*= 7250.56$.
However, no outsourcing occurs in any black swan scenario at the price of a reasonable cost $z_{ORA}^*$. See Table C.4 in the Supplementary file for more details.

\subsubsection{I1. On solving models LIP-RN \eqref{LIP-RN} and LIP-ORA \eqref{LIP-ORA} for ambiguity set $\P_{10\%}$}\label{sec:I1-RN-P10}
The experiment is similar to the one conducted for set $\P_{5\%}$.
The optimal solution cost of model LIP-RN \eqref{LIP-RN}, $z_{RN}^*=u= 7260.56$ is the same as for ambiguity $\P_{5\%}$.
Note that the setting is robust optimization and, on the other hand, $\P_{5\%} \subset \P_{10\%}$.
In fact, there are a total of eleven black swan scenarios within the ambiguity set, whose outsourcing costs show a high variability.

In the case of model LIP-ORA \eqref{LIP-ORA}, the optimal solution cost is $z_{ORA}^*= 7462.36$, which is slightly higher than for ambiguity set $\P_{5\%}$. As a result, outsourcing is completely avoided in all black swan scenarios.

\subsubsection{I1. Comparison of solutions for the ambiguity sets}\label{sec:I1-asets}
Although the number of cross-docking elements of instance I1 is relatively small, see Table \ref{I1.dim-results},
notice that the RN and ORA models' dimensions are certainly not;
as a matter of fact, they are linearly increasing with the cardinality of the ambiguity set.
Both solvers CPLEX and Gurobi obtain the optimal solution cost.
Observe the high variability on the wall time required by the solvers when dealing with the RN and ORA models for $|\P|=8$.
Notice also that CPLEX requires about 10 h 29 m to solve I1-RN-8, and
Gurobi cannot prove the optimality for instance I1-ORA-8 in the 12 h time limit, reporting a 0.03\% gap, where it is 0.04\% after 1384 s.
In any case, the ORA solutions avoid outsourcing in the black swan scenarios without a high increase in the optimal cost.

\begin{figure}[h]
\includegraphics[scale=.50]{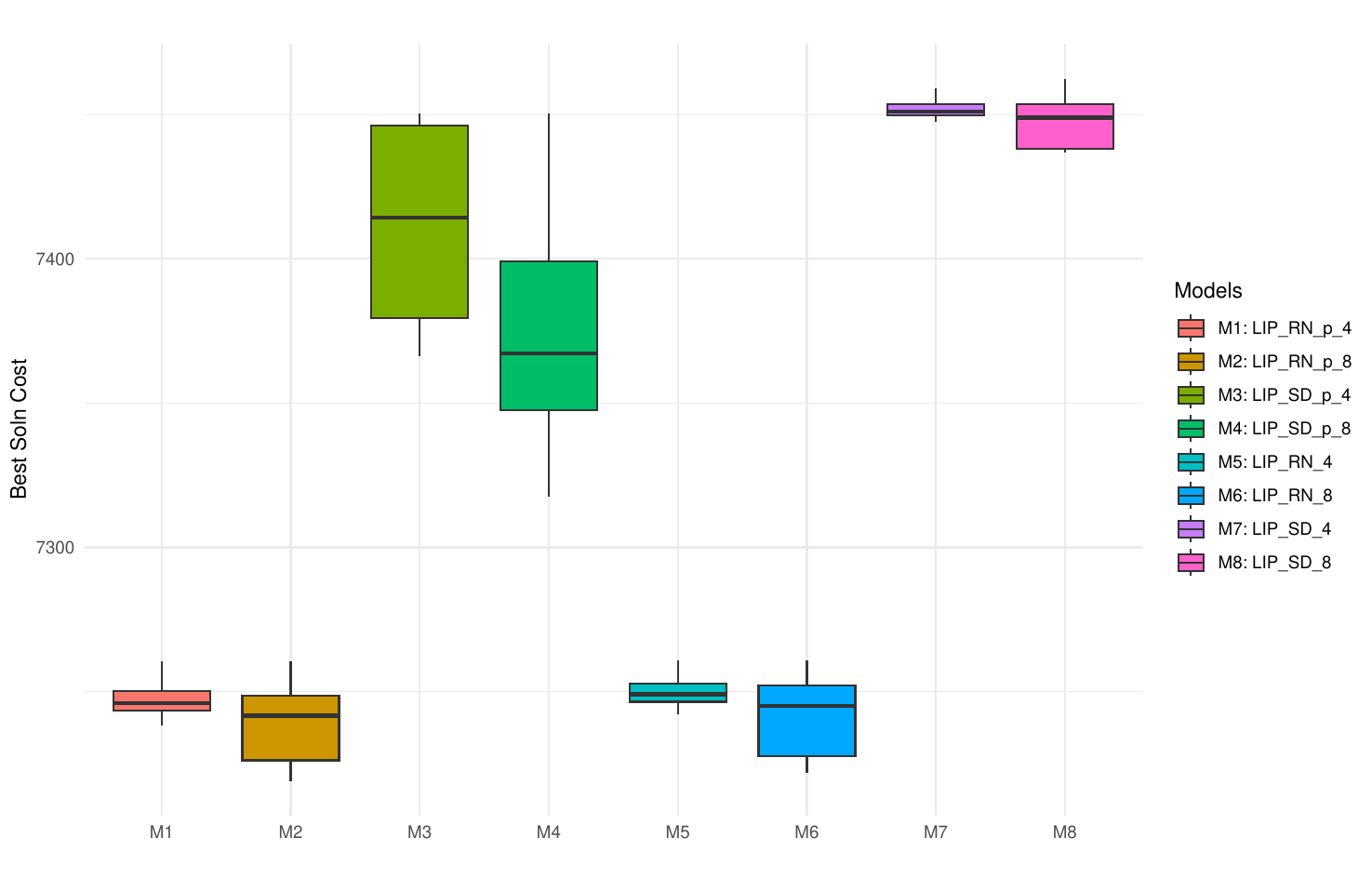}
\centering
\caption{Instance I1. DRO cost distribution retrieved from the solution of models LIP-RN \eqref{LIP-RN} and LIP-ORA \eqref{LIP-ORA} in the ambiguity sets} \label{I1-comp-solns}
\end{figure}

In addition to the comparison of computational effort, it is also interesting to analyse the CDDP cost distribution across the different members of the ambiguity set.
Fig. \ref{I1-comp-solns} depicts that cost distribution for the members of the two ambiguity sets
$\P_{5\%}$ and $\P_{10\%}$, where $|\P_{5\%}|=4$ and $|\P_{10\%}|=8$,
which is retrieved from the solution of the models LIP-RN and LIP-ORA.
The results are presented using boxplots, so that the central box in each experiment spans from the first quartile to the third one, in a interquartile mode,
the segment inside the rectangle shows the median cost.
The whiskers above and below the box show the highest and smallest costs, respectively.

The identifiers have the following meaning: M1 (M2) depicts the cost $\hat{u}^*$ of the ${\mathbb{P}}$-solution, where
model LIP-RN \eqref{LIP-RN} is independently solved for each member ${\mathbb{P}}$ in the ambiguity set $\P_{5\%}$ ($\P_{10\%}$).
The identifier M3 (M4) has a similar meaning  but replacing model LIP-RN \eqref{LIP-RN} with model LIP-ORA \eqref{LIP-ORA}.
M5 (M6) depicts the cost for each member of ambiguity set $\P_{5\%}$ ($\P_{10\%}$) that is obtained by the following two-step procedure:
First, the solution of the first stage variables is retrieved from model LIP-RN \eqref{LIP-RN}, and, second, it is fixed and a `reduced' model is independently solved for each member ${\mathbb{P}}$.
A similar scheme is used for obtaining the members' cost in M7 (M8) but replacing model LIP-RN \eqref{LIP-RN} with model LIP-ORA \eqref{LIP-ORA}.

Some observations are in order for Fig. \ref{I1-comp-solns}: The cost
that defines the upper whisker $\hat{u}^*$ is the same in M1 and M2,
and in M3 and M4, respectively, and it is higher in the solution of M3
(M4) than in M1 (M2), since the latter satisfy the ORA constraint
system \eqref{DRO-ORA}.
Fifty percent of the central costs associated with the risk-neutral
model LIP-RN reported in M5 (M6) are slightly higher than those
reported in the individual models M1 (M2); that is, the box defined by
the quartiles is slightly higher in the former case, although the
highest costs that define the upper whiskers are nearly identical.  On
the other hand, the larger the cardinality of the uncertainty set,
the greater the interquartile range or, in other words, the wider the
box in the RN models. Additionally, the cost distribution covered by the identifiers M5--M8 has the same rationale as the one for the identifiers M1--M4,
where the latter allows a lower bound on the optimal solution cost of the
original RN and ORA models. Notice, in particular, how the disruption
in cost variability fundamentally affects the LIP-ORA models, where
higher but less dispersed costs appear (M7, M8), in contrast to the
independent models solved for each member of the ambiguity set
(M3, M4). Lastly, observe that the cost distributions resulting from the ORA models are entirely above the ones obtained from the RN modeling approach.

\subsection{Instance I3}\label{sec:I3}
The instance is composed of up to 10 origins, 10 destinations, 5 strip doors and 5 stack ones.
The ambiguity sets $\P$ under consideration are $\P_{3\%}$ and $\P_{7\%}$, where $|\P_{3\%}|=3$, $|\P_{7\%}|=6$ and $\P_{3\%}\subset\P_{7\%}$;
see Appendix C.2 in the Supplementary file for their generation.

Table \ref{I3.dim-results} shows the model's dimensions of LIP-RN \eqref{LIP-RN} and LIP-ORA \eqref{LIP-ORA} as well as the main computational results and wall time required by the straightforward use of the solvers CPLEX and Gurobi.
Notice that the dimensions are approximately double than those of instance I1.
The new headings are as follows.
$\underline{z}_{L}$, LP lower bound on the optimal solution cost;
$\underline{t}_{\bullet-L}$, $\overline{z}_{\bullet}$ and $GAP_{\bullet}$,
wall time required for obtaining the LP lower bound, incumbent solution cost and optimality gap, respectively,
required by solver $\bullet$, where $\bullet\in\{cpx,grb\}$.
It is worth pointing out that the lower bound $\underline{z}_{L}$ is obtained by solving the model attached to the root node in the B\&C phase;
$\underline{t}_{cpx-L}$ and $\underline{t}_{grb-L}$ give the wall time required by CPLEX and Gurobi, respectively, for the automatically selected barrier method just before executing the crossover phase.
\begin{table}[h!]
  \begin{center}
{\scriptsize
\caption{Instance I3. CDDP models LIP-RN \eqref{LIP-RN} and LIP-ORA \eqref{LIP-ORA}: Dimensions and results of CPLEX and Gurobi}\label{I3.dim-results}
\begin{tabular}{l|cccc|c|ccc|ccc}
\hline
  $Inst$-$|\P|$ & $m$ & $n01$ & $nc$ & $nz$
                   & $\underline{z}_{L}$
                   & $\underline{t}_{cpx-L}$ & $\overline{z}_{cpx}$ & $GAP_{cpx}$
                   & $\underline{t}_{grb-L}$ & $\overline{z}_{grb}$ & $GAP_{grb}$\\
\hline
I3-RN-3& 28995 &3110 &78031 &230823 &7984.59 &3 &8451.89 &5.5 &20&8339.79  &4.3\\
I3-ORA-3& 29131 &3113 &78094 &232813 &8198.65 &2 &8526.87  &3.8 &6&8503.66  &3.6\\
I3-RN-6& 57978 &6170 &156061 &461546 &7991.35 &7 &8597.93 &7.0 &229&8315.29  &3.9\\
I3-ORA-6& 58249 &6176 &156187 &465526 &8230.65 &9 &8695.13 &5.3 &23&8512.42  &3.3\\
\hline
\end{tabular}
}
\end{center}
\end{table}

The first observation about the elements in Table \ref{I3.dim-results} is that the models are very difficult to solve;
as a matter of fact, they cannot be solved up to optimality in the given time limit.
Notice the small time required by both solvers for obtaining the LP bound $\underline{z}_{L}$ versus the 12 h time limit that is reached while solving the other nodes in the B\&C phase.
It is also worth observing that both solvers provide the \textit{same} LP bound in spite of their potentially different preprocessing approaches.
The optimality $GAP$ is expressed as $100.\frac{\overline{z}-\underline{z}}{\overline{z}}$,
as typically done in the commercial solvers,
where $\overline{z}$ is the incumbent cost and $\underline{z}$ is the smallest cost of the solutions retrieved from the submodels in the set of active B\&C nodes at the optimization's interruption time.
It can be observed that $GAP_{grb}$ is smaller than $GAP_{cpx}$.
Additionally, notice that the ORA-incumbent cost is higher than the RN-incumbent one.

Observe that there are ten black swan scenarios with high costs in the ambiguity set $\P_{3\%}$ for  LIP-RN,
due to outsourcing, and a very small weight in the three ambiguity set members (see Table C.8 in Appendix C).
These outsourcing costs are prevented in the LIP-ORA model by introducing an appropriate modeler-driven pair of parameters' triplex
$(\iota,\overline{s},\rho)$. As a result, the black swan scenarios do not have negative implications for the ORA solution obtained by both solvers,
although the cost increases moderately (see Tables 7 and C.9).
The rationale for set $\P_{7\%}$ is similar to the one presented above for $\P_{3\%}$.

\begin{figure}[h]
\includegraphics[scale=0.50]{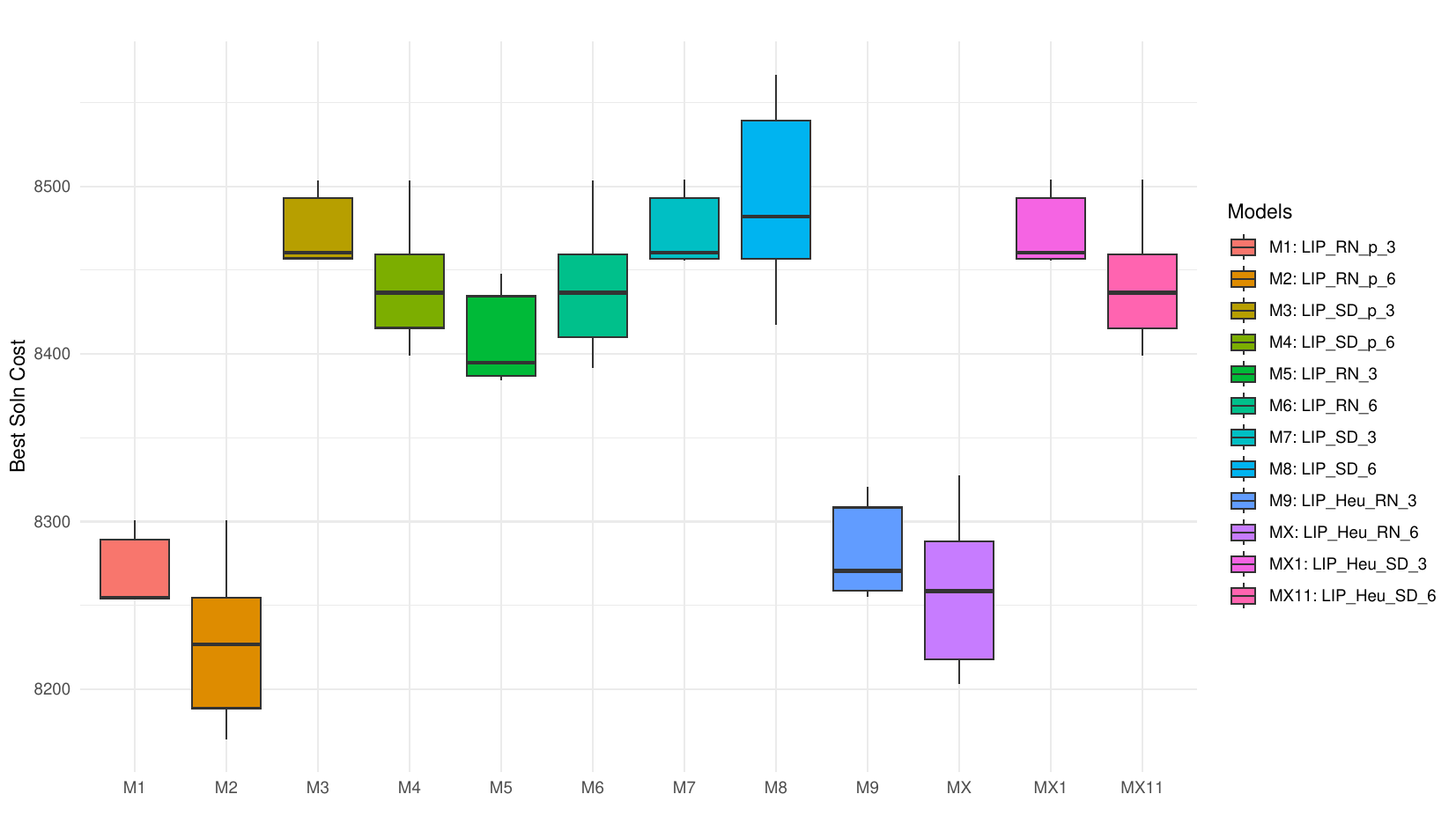}
\centering
\caption{Instance I3. DRO cost distribution retrieved from the solution of models LIP-RN \eqref{LIP-RN} and LIP-ORA \eqref{LIP-ORA} in the ambiguity sets} \label{I3-comp-solns}
\end{figure}

Fig. \ref{I3-comp-solns} depicts a comparison of the CDDP cost distribution of the ambiguity set members with the same meaning as Fig. \ref{I1-comp-solns}
for M1--M8, where now the ambiguity sets are $\P_{3\%}$ and $\P_{7\%}$, $|\P_{3\%}|=3$ and $|\P_{7\%}|=6$.
The group of four identifiers, M9--MX and MXI–MXII, serves a similar purpose to the group consisting of M5–M6 and M7–M8. The key difference is that in the former group, the proposed matheuristic scheme replaces the direct use of CPLEX for solving the corresponding models, see Section \ref{sec:H}.
Notice that the experiment has the same rationale for M5--M8 and
 M9--MXII. A comparison between the solutions obtained using CPLEX
  and those derived from the matheuristic reveals a
 substantial difference in the cost distributions for the RN-related
 identifiers, M5 (M6) and M9 (MX), as well as for the ORA-related
 ones, M7 (M8) and MXI (MXII). It highlights the poor performance of
 the CPLEX solver in obtaining high-quality solutions. With regard to
 the solutions of the RN versions, it is observed that the cost
 defining the upper whisker $\hat{u}^*$ is identical in M1 and M2, and
 lower than those obtained in M9 and MX. In the latter two models, the
 costs differ, since the optimization has been performed across all
 elements of the uncertainty sets with different cardinalities.
 On the other hand, the cost distribution obtained from the solutions
 provided by the matheuristic for the model LIP-ORA \eqref{LIP-ORA},
 MXI (MXII), practically replicates the one obtained individually in
 M3(M4), with these costs being significantly higher than those obtained in the RN models M1 (M3) and M9(MX).

Finally, let us point out an analysis of the reliability of the
ambiguity sets $\P_{3\%}$ and $\P_{7\%}$, where $\P_{3\%} \subset
\P_{7\%}$, by using  the information provided by the identifiers M5
and M6 (or, M9 and MX).
Observe that the feasible set in the former is within the one in the latter, the median cost is higher in the latter as well as it is the whisker above the respective box.
It gives some hints on the higher accuracy of $\P_{7\%}$ versus $\P_{3\%}$,
i.e., it looks as the highest cost in the former is more accurate than the one in the latter;
see also the results shown in Table  \ref{I3.dim-results}.
A similar analysis can be done for the identifiers M7 and M8 (or MXI and MXII).

\subsection{Instance I7}\label{sec:I7}

The instance is composed of up to 20 origins, 20 destinations, 10 strip doors and 10 stack ones.
The ambiguity sets $\P$ under consideration are $\P_{3\%}$ and $\P_{7\%}$, where $|\P_{3\%}|=3$, $|\P_{7\%}|=6$ and $\P_{3\%}\subset\P_{7\%}$;
see Appendix C.3 for their generation.

\subsubsection{I7. On solving models LIP-RN \eqref{LIP-RN} and LIP-ORA \eqref{LIP-ORA}}\label{sec:I7-RN-P3}
Table \ref{I7.dim-results} shows the model dimensions for LIP-RN \eqref{LIP-RN} and LIP-ORA \eqref{LIP-ORA} as well as the main computational results.
Notice that the dimensions are one order of magnitude higher than the ones of instance I3.
So, the models are much more difficult to solve and, so, the optimization also reaches the time limit (12 h), if any.
It can be observed that the CPLEX solution for the RN version is very bad, it has been obtained in the solver's preprocessing;
it is worthless in this respect.
On the other hand, the solution provided by Gurobi for the RN  model is an acceptable one:
observe the small time required by the barrier method for obtaining the LP cost without using the crossover scheme,
since otherwise the time would be one order of magnitude higher.
It is worth pointing out that none of the two state-of-the-art solvers could find any feasible solution in the 12 h wall time limit for the ORA models.
CPLEX is even unable to provide an LP bound for the largest instance I7-ORA-6.
It can be deduced from the results shown in Tables \ref{I3.dim-results} and \ref{I7.dim-results} that Gurobi 11.0.3 outperforms CPLEX 22.1.1 whenever a feasible solution is found.
\begin{table}[h!]
  \begin{center}
{\scriptsize
\caption{Instance I7. CDDP models LIP-RN \eqref{LIP-RN} and LIP-ORA \eqref{LIP-ORA}: Dimensions, and CPLEX and Gurobi results\\
Note: $GAP_{cpx}=99.94$ for I7-RN-3 and I7-RN-6; $\,\,\, -$, No feas soln found}\label{I7.dim-results}
\begin{tabular}{l|cccc|c|cc|ccc}
\hline
  $Inst$-$|\P|$ & $m$ & $n01$ & $nc$ & $nz$
                   & $\underline{z}_{L}$
                   & $\underline{t}_{cpx-L}$ & $\overline{z}_{cpx}$
                   & $\underline{t}_{grb-L}$ & $\overline{z}_{grb}$ & $GAP_{grb}$\\
\hline
I7-RN-3 &210865 &12970 &969436  &2568280  &20723.10 &  304 & 32419500 & 220 & 23394.69 &
11.4\\
I7-ORA-3 & 211121&12973 &969559  &2575214 &20947.13 & 26844 &      $-$  & 130 &      $-$ & $-$\\
I7-RN-6 &421708 &25840 &1938871 &5136361  &21612.95 &  233 & 41457800  & 205 & 24224.44 & 10.8\\
I7-ORA-6 &422219 &25486 &1939479 &5150231  &21612.95 & $-$  &   $-$ & 244 &      $-$ & $-$\\
\hline
\end{tabular}
}
\end{center}
\end{table}
\subsection{Performance of the matheuristic scheme}\label{sec:H}
The section presents the performance of the matheuristic scheme introduced in Sections \ref{sec:rn-lb-ub} and \ref{sec:ORA-lb-ub},
while solving the CDDP models LIP-RN and LIP-ORA, respectively.
For this purpose, the difficult instances I3 and I7 are considered, which model dimensions are shown in Tables \ref{I3.dim-results} and \ref{I7.dim-results}, respectively.

Table \ref{heur-results} shows the main results obtained by the matheuristic.
Its headings are as follows:
$\underline{z}_L$, LP lower bound on the optimal solution cost taken from Tables \ref{I3.dim-results} and \ref{I7.dim-results};
$\underline{z}_H$, lower bound obtained by the matheuristic;
$\underline{t}_L$ and  $\underline{t}_H$, related wall times;
$\hat{C}_{1,H}$ and $\hat{F}_H$, first stage cost and second stage expected cost, respectively
(they compose the total cost $\overline{z}_{H}$ of the incumbent solution for the ambiguity set member provider);
$\overline{t}_{H}$, wall time required for obtaining $\overline{z}_{H}$, where time $\underline{t}_H$ has been added;
$GAP_H\%$, optimality gap of the incumbent cost to be expressed as $100\cdot{{\overline{z}_H-\underline{z}'}\over {\overline{z}_H}}$,
in the spirit of the solvers, where $\underline{z}' = max \{ \underline{z}_L, \, \underline{z}_H \}$; and
$GR_H$, goodness ratio of cost $\overline{z}_H$ with respect to the Gurobi incumbent cost $\overline{z}_{grb}$.
It is worth pointing out that
$\underline{z}'=\underline{z}_H$ for I3-ORA-3 and I3-RN-6, otherwise,
$\underline{z}'=\underline{z}_L$.

It is remarkable the usefulness of the first stage alternative solutions  for obtaining the incumbent solution cost $\overline{z}_H$,
given the wall time, optimality gap and goodness ratio as shown in Table \ref{heur-results}.
Notice that the first stage solutions are retrieved from the cluster decomposition-based lower bound approaches, see Sections \ref{sec:rn-lb-ub} for RN and \ref{sec:ORA-lb-ub} for ORA.

\begin{table}[h!]
\begin{center}
{\small
\caption{Instances I3 and I7. CDDP models LIP-RN \eqref{LIP-RN} and LIP-ORA \eqref{LIP-ORA}: Matheuristic results}\label{heur-results}
\begin{tabular}{l|cc|cc|cc|cc|cc}
\hline
$Inst$-$|P|$ & $\underline{z}_L$ & $\underline{t}_L$ & $\underline{z}_H$ & $\underline{t}_H$ &  $\overline{z}_H$ & $\overline{t}_H$
& $\hat{C}_{1,H}$ & $\hat{F}_H$& $GAP_H$ & $GR_H$\\
  \hline
I3-RN-3 &  7984.59 &   20 &6456.43 &32&  8320.96 &   161 &1849&6471.96& 4.04 & 0.998 \\
I3-ORA-3 &  8198.65 &  6 & 8226.63  & 360&  8503.66 &   402  &2082&6421.66& 3.26 & 1.000 \\
I3-RN-6 &  7991.35& 229 &8130.70 &  234&  8327.80 &  564 &1849&6478.80&  2.37 & 1.001 \\
I3-ORA-6 &  8230.65 &  220 &8226.63 &636&  8503.94 &  793  &2089&6414.94& 3.21 & 0.999 \\
I7-RN-3 & 20723.10 &  220 &17418.44 &14400& 22402.18 &  14787 &7157&15245.18&  7.49 & 0.957 \\
I7-ORA-3 & 20947.13 &  123 &19435.35 &14400& 22402.69 &14876   &7284&15118.69&  6.49 & $\infty$ \\
I7-RN-6 & 21612.95 &  205 &20078.54 & 14400  &23700.67&14875&7784&15916.67&8.81 & 0.971 \\
I7-ORA-6 & 21612.95 &  206 &20729.06 &14400& 23939.31 & 14845     &8023&15916.31&   9.71   & $\infty$ \\
\hline
\end{tabular}
}
\end{center}
\end{table}

Table \ref{heur-results} shows similar types of results of the matheuristic for instance I3 as those shown in Table \ref{I3.dim-results} while solving it by straightforward use of the two commercial solvers.
Additionally, similar observations can be made on the results
retrieved from the solutions of models LIP-RN \eqref{LIP-RN} and
LIP-ORA \eqref{LIP-ORA} for instance I7, see Tables \ref{I7.dim-results} and \ref{heur-results}.
Notice that the black swan scenarios have outsourcing in the RN solutions for the four instances I3-RN-3, I3-RN-6, I7-RN-3 and I7-RN-6.
On the other hand, it is also remarkable that the ORA solutions have no outsourcing.
Additional results are shown in Tables C.9, C.10 and C.11 for Instance I3 and Tables C.15 and C.16 for Instance I7.

\textbf{Parallelism.} One should notice the embarrassing parallelism on obtaining the bounds provided by the matheuristic scheme.
Notice that the lower bounds \eqref{lb-RN} for RN and \eqref{lb-ORA} for ORA are obtained by independently solving and, then, in parallel,
the $|\C|$ $c$-submodels \eqref{RN-SCD-c} and $c$-submodels \eqref{ORA-SCD-c}, respectively, plus the LP relaxation of the related models \eqref{LIP-RN} and \eqref{LIP-ORA}.
So, $\underline{t}_{H}$ is the highest wall time for solving the most demanding related $c$-submodel.
In a similar way, the upper bounds \eqref{ub-RN} for RN  and \eqref{ub-ORA} for ORA are obtained by independently solving and, then, in parallel,
the  $\omega$-submodels \eqref{RN-cdap} and \eqref{ORA-cdap},
according to the available hardware structure.
Those submodels result from decomposing the models \eqref{LIP-RN} and \eqref{LIP-ORA}, respectively, based on the appropriate fixing of the first stage variables.
So, $\overline{t}_{H}$ is the highest wall time for solving the most demanding related submodel for each scenario in the matheuristic versions, where $\underline{t}_{H}$ is added.
Note: The time limit for solving each of the submodels has been set to 4 h.

\section{Conclusions and future research agenda}\label{sec:con}
In this work, we have introduced a two-stage distributionally robust
optimization approach (CDDP-DRO) to address the cross-dock platform design problem
under uncertainty. Its risk-neutral version considers uncertain freight
flows and aims to minimize expected handling costs while ensuring
robust performance across a range of plausible probability
distributions.
Additionally, an outsourcing risk-averse measure has been proposed to prevent the negative effects of black swan scenarios, defined here as the occurrence of `outsourcing' solutions.

All of the elements considered in the definition of CDDP-DRO
contribute to the overall complexity of the problem, which is already
NP-hard with a quadratic objective function. To the best of our knowledge, this challenging problem has not been previously addressed in the literature, despite its broad applicability across various fields.
We have proposed a \textit{min-max}-based scheme that
decomposes each element of the ambiguity set into scenario clusters,
enabling the derivation of a robust solution.

Numerical experiments are reported for comparing the solutions provided by
the straightforward use of the state-of-the-art solvers CPLEX and
Gurobi and validating the proposed matheuristic scheme.
 The testbed in the computational experiment was composed of extended
versions of the instances I1, I3 and I7, originally presented in
\cite{cddp-ts24}. The model dimensions can reach up to around 420,000 constraints, 25,000 binary variables, and 1,900,000 continuous variables, with complex combinatorial structures.

A computational comparison between the straightforward use of both
commercial solvers was performed by using the RN and ORA versions of the model.
For medium-sized instances, both solvers reach the 12 h wall time limit without proving the incumbent DRO solution's optimality;
in fact, the optimality gaps for CPLEX and Gurobi are up to 7.0\% and 4.3\%, respectively.
For large-sized RN instances, CPLEX fails to deliver meaningful
results, while Gurobi yields quasi-optimal solutions with gaps up to 11.4\%.
Moreover, neither of these solvers is able to provide a feasible solution for the large-sized ORA instances.

On the other hand, the proposed matheuristic scheme provides an incumbent solution for all tested runs in a very reasonable wall time;
its optimality GAP ranges from 2.37\% to 9.71\% for medium- and
large-sized instances. Interestingly, the goodness ratio of the cost provided by the matheuristic solution versus that of Gurobi ranges from 0.957 to 1.001
for the runs where the latter provides a solution.

Finally, our future research agenda will include the following items:
(a) a generalization of the proposed matheuristic scheme for solving general two-stage DRO MILP problems, and then benchmarking against approaches such as the one presented in \cite{Jin26},  and
(b) development of a multistage multi-horizon DRO MBQ approach for CDDP designing.

\bigskip
\noindent\textbf{CRedIT author statement}

\textbf{Laureano F. Escudero}: Conceptualization, Methodology, Investigation, Writing original draft preparation, Writing-reviewing and editing, Supervision.
\textbf{M. Araceli Gar\'in}: Methodology, Formal\break analysis, Investigation, Data curator, Software, Writing-reviewing and editing.
\textbf{Aitziber Unzueta}: Methodology, Investigation, Data curator, Software, Writing-reviewing and editing, Validation.

\bigskip
\noindent\textbf{Acknowledgments}

The authors wish to express their thanks to Antonio Alonso-Ayuso, Juan Francisco Monge, Domingo Morales and Leandro Pardo for fruitful discussions on the generation of the ambiguity set.

This research has been partially supported by the projects
RTI2018-094269-B-I00 and PID2021-122640OB-I00 (L.F. Escudero), Grant PID2023-147410NB-100 funded by MCIN/AEI/10.13039/ 501100011033
and by ERDF/EU and Grupos de Investigaci\'on EOPT (IT-1494-22) and
MOMASA (IT-2072-26) from the Basque Government
(M.A. Gar\'in and A. Unzueta).

\bigskip
\noindent\textbf{Declarations}

The authors have no competing interests to declare that are relevant to the content of this article. The data that support the findings of this study can be available upon reasonable request to the authors.

\bigskip




\newpage
\centerline{\bf \large Cross-Dock Door Design under uncertainty. A
  two-stage DRO-based}

\centerline{\bf \large lower- and   upper-bounding  scheme. Supplementary file: Appendices}





\bigskip
\normalsize

\goodbreak
\appendix{\bf\normalsize}

\setcounter{table}{0}
\section{An overview of applications of Distributional Robust Optimization}
Some applications of Distributional Robust Optimization  (DRO) are as follows:
\cite{SaifDelage21} presents a DRO approach for capacitated facility location problem.
\cite{Arrigo22} present a multistage distributionally robust chance-constrained mixed integer linear programming approach with a Wasserstein ambiguity set, for a mix of conventional and renewable energy generators system design.
\cite{Hosseini-Nodeh22,Mei22,Zhu09} present approaches for portfolio management.
\cite{Jiang24} present an approach for the distribution of potassium
fertilizer from suppliers to customers through a set of cross-docks along a
time horizon under uncertain demand.
The deterministic LP model is converted into a second-order cone program by using a DRO scheme to deal with the uncertainty.
Another DRO approach is presented in \cite{Lu24} for parallel identical machines scheduling where the job processing time is uncertain.
A mean--mean absolute deviation-based ambiguity set is considered.
A DRO support vector clustering LP model for supply chain management is studied in \cite{Sun25};
the uncertainty relies on the product's demand,
where the true probability distribution is assumed to belong to an ambiguity set to be estimated by considering a demand data set.
In particular, onmichannel strategies are studied when selling the product through online and physical stores.
A multistage DRO is presented in \cite{Yang24} for a medical supplies distribution network with uncertain demand.
An ambiguity set is considered, being based on the mean absolute deviation and an autoregressive approach integrating a moving average method.
The decisions at each node in the tree only depend on the parameter's realizations at the node and its immediate ancestor.
A brief overview of additional DRO applications follows:
biomass technology with agricultural waste-to-energy network design in \cite{NingYou19},
machine scheduling and humanitarian logistics in \cite{Noyan22},
multistage water allocation from rivers to consumption areas in \cite{ParkBayraksan23},
wind farm layout in \cite{AgraCervera24},
berth allocation in \cite{AgraRorigues24},
disaster relief planning in \cite{ElTonbary24} and \cite{Wang23},
shortest path problem in \cite{Filippi25},
closed-loop supply network design for recycling and manufacturing in \cite{Gao24} and \cite{Leng25},
online portfolio selection for an environmental, social and governmental system in \cite{Guo25},
surgery scheduling in \cite{Li25},
principal--agent contract design in \cite{Liu25},
wind power-type energy systems management in \cite{Ou26} and its literature review,
locating-routing in urban search and rescue operations in \cite{Sarmadi25},
 ranking priorities of social networks in \cite{Shao25},
construction and demolition waste disposal facilities in \cite{Xin25},
ensemble forecasting modeling in \cite{Yuan25},
train make-up and passenger flow control on an overcrowded metro line in \cite{Yang26},
mobile energy generators capacity and storage planning in \cite{Zhang25},
production pricing, quality and quantity design and adjustments in \cite{ZhangTian24},
and others. See also the recent review in \cite{Lin22}.
A framework for modeling and approximately solving DRO problems with decision-dependent information discovery is proposed in \cite{Jin26}.

\setcounter{table}{0}
\section{A procedure for generating the Wasserstein ambiguity set}\label{sec:aset}
This appendix presents the features of the procedure that is proposed for generating the ambiguity set for the second stage uncertain parameters in the experiment.

The main assumptions considered on generating the ambiguity set are as follows:
\begin{itemize}\parskip 0.5mm
\item Assumption 1: The set of uncertain parameters is modeled as a set of independent random variables.
\item Assumption 2: There is no information on the true probabilistic distribution that is followed by the uncertainty.
However, empirical information on their nominal distribution (ND) is available. The uncertain realizations are discretized into a set of scenarios denoted by $\Omega_{\mathbb{P}_N}$ and the scenario weights are also available.
\end{itemize}

Roughly speaking,  the process for generating the ambiguity set of each uncertain parameter has the two following steps:
\begin{enumerate}\parskip 0.5mm
\item Generate new realizations and weights for the uncertain parameters. \\
To this purpose,   the mean ($\mu_N$) and variance (${\sigma}^2_N$) of the realizations of each parameter under  nominal distribution $\mathbb{P}_{N}$  are first computed.   Subsequently, the realizations of the uncertain parameters from the nominal distribution under scenario $\omega \in \Omega_{\mathbb{P}_N}$, denoted by $\xi^{\omega}_N$, are projected using the cumulative distribution function of each available probability distribution $\mathbb{Q} \in \Q$. These projections are obtained as follows,
\begin{equation}
u^{\omega}_{\mathbb{Q}}=F_\mathbb{Q}(\xi^{\omega}_N;\mu_N,{\sigma}^2_N) \quad  \forall \omega \in \Omega_{\mathbb{Q}}, \quad  \mathbb{Q} \in \Q
\end{equation}

The next step consists of generating perturbations. To this end, $\epsilon$-perturbations are sampled from a Normal Distribution, $\epsilon \sim N(0,\sigma_\epsilon^2)$, and applied to obtain perturbed values $u'^{\omega}_\mathbb{Q} = u_\mathbb{Q} + \epsilon$. The resulting values of $u_\mathbb{Q}$ are then truncated to maintain admissible bound. $u'^{\omega}_\mathbb{Q}$ is set to zero whenever a non-positive value is obtained and set to one whenever the value exceeds one.
Finally, the realization of each uncertain parameter under each scenario is generated by applying the inverse cumulative probabilistic function for each available probability distribution $\mathbb{Q} \in \Q$ can be expressed as
\begin{equation}
\xi^{\omega}_\mathbb{Q} =F_\mathbb{Q} ^{-1}(u'^{\omega}_\mathbb{Q} ;\mu_N,{\sigma}^2_N) \quad \forall \omega \in \Omega_{\mathbb{Q}}
\end{equation}

Given Assumption 1, that is, the independence of the uncertain parameters (i.e., random variables)
the likelihood function associated with each scenario, say $L^{\omega}$, together with the scenario’s weight under the nominal distribution $w^\omega_N$,  is used to derive the corresponding normalized scenario weights $w^{\omega}$ under distribution $\mathbb{Q}$, satisfying $\displaystyle \sum_{\omega\in\Omega_\mathbb{Q}}w^\omega=1.$  This reweighting of scenarios—typical of entropic or Kullback-Leibler (KL)-type divergences—means that the method is no longer purely Wasserstein-based, where the formulation involves only generating new scenarios while keeping the same weights.

This combination enables richer modeling of scenarios by coupling the geometric robustness of Wasserstein with the statistical robustness of KL divergence. It precisely combines adversarial training (Wasserstein) and importance weighting (KL), and we find it particularly suitable for analyzing risk aversion in our setting, where worst-case scenarios not only receive greater weight but also become more adverse.

\item Select the Wassertein ambiguity set $\P$ using previously generated elements.\\
The aim is to generate the ambiguity set $\P$, so that the realization of the uncertain parameters is represented in a more accurate way. For that purpose,
we select those probability distributions $\mathbb{Q}$ which realizations are closest to those obtained from the Nominal distribution, $\mathbb{P}_N$.  Wassertein distance is employed as a selection criterion, meaning that the weighted distance,  \textit{proximity} $l_{\mathbb{Q}}^r$, between
$( \xi_\mathbb{Q}^\omega, w^\omega \, \, \forall \omega\in\Omega_\mathbb{Q} )$ in perturbation $\mathbb{Q}$ and
$( {\xi}_{N}^{\omega'}, w^{\omega'} \, \, \forall \omega'\in\Omega_{\mathbb{P}_N})$ in nominal distribution  $\mathbb{P}_N$ is small enough
(i.e., no greater that, say, $\theta$).
 The minimization problem for the selection of the desired ambiguity set  $\P$ can be expressed as
\begin{equation}\label{wasser}
\begin{array}{lll}
&& \displaystyle l_{\mathbb{Q}}^r =
min_{\eta^{\omega,\omega'} \geq 0} \sum_{\omega\in\Omega_\mathbb{Q}} \sum_{\omega'\in\Omega_{\mathbb{P}_N}} d^r(\xi_\mathbb{Q}^\omega,{\xi}_{N}^{\omega'})\eta^{\omega,\omega'}\\
&& \text{ s.to } \displaystyle \sum_{\omega'\in\Omega_{\mathbb{P}_N}}\eta^{\omega,\omega'}  =w^\omega           \quad \forall \omega\in\Omega_\mathbb{Q}\\
&&               \displaystyle \sum_{\omega\in\Omega_\mathbb{Q}} \eta^{\omega,\omega'} =w^{\omega'} \quad \forall \omega'\in\Omega_{\mathbb{P}_N},
\end{array}
\end{equation}
where $\eta^{\omega,\omega'}$ is the `transportation' variable from scenario $\omega$ to scenario $\omega'$ and
$d^r(\xi_\mathbb{Q}^\omega,{\xi}_{N}^{\omega'})$ is the distance between perturbation $\mathbb{Q}$ and nominal distribution $\mathbb{P}_N$, and $r\in \{1,2, \infty\}$. Let $\P$ be the subset of the candidate members such that $l_\mathbb{Q}^r \leq \theta$.

\end{enumerate}

Before presenting the detailed generation of the ambiguity set for the instances in the experiment, note that the parameters that represent the inbound volume entering the platform and the outbound volume exiting from, say $S_m^\omega$ and $R_n^\omega$, respectively, and the commodity handling and transportation cost, say $G_{minj}^\omega$, are directly obtained from $H_{mn}^\omega$, the commodity volume to be consolidated in the platform from origin node $m$ to destination node $n$. Then, the uncertain parameters to which the proposed methodology has been applied are: {\it a)} the commodity volume  $H_{mn}^\omega$,  and  {\it b)} the  capacity disruption fraction of the strip and stack doors $D_i^\omega$ and $D_j^\omega$, respectively. Taking into account the stochasticity of these second stage parameters, different probability distributions were used to generate the ambiguity set. On the one hand, the Normal (1), Weibull (2), Gamma (3) and Lognormal (4) distributions,  denoted by $\mathbb{Q}=q\in\{1,2,3,4\}$,  were used to generate the different perturbations of the realizations $H_{mn}^\omega$. On the other hand, the perturbations of the  capacity disruptions $D_i^\omega$ and $D_j^\omega$ were generated using the Uniform distribution.

\setcounter{table}{0}
\section{Ambiguity sets generation and DRO solutions for instances I1, I3 and I7}\label{sec:aseti1i3i7}
The appendix details the procedure for the generation of the ambiguity sets for the instances in the experiment; it is based on the procedure presented in \ref{sec:aset}.
The statistical handling  of the data throughout this section was carried out
using the open source R statistical package, see \cite{r}.

\subsection{Ambiguity sets generation for instance I1}\label{sec:dataI1}
The smallest instance in the experiment is composed of up to 8 origins and 8 destinations, and 4 strip doors and 4 stack doors.
In particular, the ND is defined by a set of $|\Omega|=5$ equiprobable scenarios in $|\L|=1$ group of scenarios,
and the four strip and stack doors have no capacity disruption (i.e., $D_i^\omega=D_j^\omega=0 \, \forall \omega\in\Omega$).
The related non-zero realizations (i.e., commodity volumes $H_{mn}^\omega$) are chosen \textit{at random},
for scenario $\omega\in\Omega$ in ND.

Thus, a candidate ambiguity set $\P'_q$ is built for each distribution $q \in \Q$ and the weight $w^\omega$ associated with each scenario in each perturbation ${\mathbb{P}} \in \P'_q$ are computed. Once those elements are obtained as well as the distance matrices between the ND realization and the realizations of the uncertain parameters, the `transportation' model (B.3) is solved using the  \textit{proximity} $\ell_{\mathbb{P}}^r$, for $r=2$. On the generation of this ambiguity set, the parameter $\sigma_ \epsilon=0.05$ (standard deviation in the normal distribution of the $\epsilon$-perturbation), and the cardinalities $|\Q|=4$ and $|\P'_q|=20$ are considered, so $|\P'_q| \times |\Q|=80$ `transportation' models are solved.

Then, a set of perturbations ${\mathbb{P}}$ is to be selected.
For this purpose, the modeler-driven parameter, the \textit{radius} $\theta$, is determined based on specific percentiles from the proximity set.
Table \ref{Tab2.I1} shows the statistics of $\ell_{\mathbb{P}}^2$ for the candidate members to be included in the ambiguity set $\P$.
The set $\P=\P_{10\%}$ corresponds to the first decile $\theta_{10\%}=9.338$ as the \textit{radius};
it is composed of the eight consecutively numbered perturbations $\mathbb{P}\in\{1,2,...,8\}$, and, in this case, all of them are Normal.
On the other hand, the ambiguity set $\P_{5\%}$ does it for the first 5-centile reduced \textit{radius} $\theta_{5\%}=9.145$;
it is composed of the first four perturbations $\mathbb{P}\in\{1,2,3,4\}$, so that $\P_{5\%}\subset\P_{10\%}$.

\begin{table}[!h]
\begin{center}
{\small
\caption{Instance I1. Statistics for the candidates for being the ambiguity set}\label{Tab2.I1}
\begin{tabular}{cccccccc}
  \hline
 Min&1st 5cent&1st 10cent&1st Qu&Median& Mean&3rd Qu&Max\\
  \hline
8.855&9.145&9.338&10.048&10.411&10.545&11.037&12.570\\
  \hline
 \end{tabular}
 }
\end{center}
\end{table}
\subsubsection{Instance I1. More details on generating ambiguity sets $\P_{5\%}$ and $\P_{10\%}$ and DRO solutions}\label{sec:I1details}
Table \ref{Tab1.I1} summarizes some other statistics for $\ell_{\mathbb{P}}^2$ by considering the four distribution assumptions in set $\Q$.
The statistics are the minimum, first quartile, median, mean, third quartile and the maximun value in the frequency distribution of the proximities.
It can be seen that the perturbations in the Normal distribution are the nearest ones to ND.
\begin{table}[!h]
\begin{center}
{\small
\caption{I1. Statistics for \textit{proximity} $\ell_{\mathbb{P}}^r \, \forall \mathbb{P}\in\P'_q, \, q\in\Q$}\label{Tab1.I1}
\begin{tabular}{lrrrrrr}
  \hline
\text{Distribution}&  Min & 1st Qu & Median& Mean & 3rd Qu & Max\\
   \hline
  Normal&8.855&9.178&9.367&9.407&9.529&10.210\\
  Weibull&11.540&11.830&11.930&11.970&12.050&12.570\\
  Gamma&9.936&10.267&10.411&10.401&10.560&10.775\\
  Lognormal&10.040&10.270&10.420&10.400&10.520&10.870\\
  \hline
\end{tabular}
}
\end{center}
\end{table}

Tables \ref{Tab3.I1.RN.P5} and \ref{Tab3.I1.ORA.P5} show the weights
and optimal solution of models LIP-RN (8) and LIP-ORA (19),
respectively, for  ambiguity set $\P_{5\%}$.
Notice in the optimal solution of model LIP-RN that scenario $\omega=5$ is a black swan in all ambiguity set members and
it is scenario $\omega=1$ in members $\mathbb{P}=2$ and $3$.
The  high cost is due to outsourcing; they also share a very small weight.
For the ORA parameters' triplex ($\iota=7600$, $\overline{s}=2500$, $\rho=2000$),
the optimal solution cost $z_{ORA}^*= 7459.31$ of model LIP-ORA (19) for ambiguity set $\P_{5\%}$ is higher than the RN counterpart $z_{RN}^*=7260.56$.
It is composed of the first stage infrastructure building cost $\hat{C}_1=1764$  and
the second stage highest expected cost $5695.31$ (where the selected ambiguity set member is $\mathbb{P}=4$, since $\gamma_{\mathbb{P}}=1$).
It can be observed that the black swan scenarios have not any more so negative implication (i.e., outsourcing, in this case)
in the solution of model LIP-ORA (19) at a reasonable cost $z_{ORA}^*$.
\begin{table}[!h]
\begin{footnotesize}
\begin{center}
\caption{I1. Weights and optimal solution of model LIP-RN (8) for the cross docking node$-$door assigned costs in scenario $\omega$, for $\omega\in \Omega_{\mathbb{P}}, \, \mathbb{P}\in\P_{5\%}$,}
\label{Tab3.I1.RN.P5}
\begin{tabular}{|cc|cc|}
    \multicolumn{2}{c}{$\mathbb{P}=1$}&\multicolumn{2}{c}{$\mathbb{P}=2$}\\
  \hline
 $w^1=5.841132e-06 $ &$\hat{F}^1=$5388&$w^1=3.909180e-05$&$\hat{F}^1=$12401188\\%
  $w^2=6.082577e-07$ &$\hat{F}^2=$5555& $w^2=4.025044e-02$&$\hat{F}^2=$5450\\
 $w^3=7.230249e-01$  &$\hat{F}^3=$5639& $w^3=7.701221e-01$&$\hat{F}^3=$5639\\
  $w^4=2.161435e-01$ &$\hat{F}^4=$6033& $w^4=1.895883e-01$&$\hat{F}^4=$6041\\
  $w^5=2.404819e-12$ &$\hat{F}^5=$14000416& $w^5=9.389944e-25$&$\hat{F}^5=$14400000\\
    \hline
    \multicolumn{2}{c}{$\mathbb{P}=3$}&\multicolumn{2}{c}{$\mathbb{P}=4$}\\
      \hline
 $w^1=3.796548e-07$&$\hat{F}^1=$14600000& $w^1=7.957016e-06$&$\hat{F}^1=$5380\\
 $w^2=3.658705e-02$&$\hat{F}^2=$5483& $w^2=5.807343e-02$&$\hat{F}^2=$5451\\
 $w^3=7.54214e-01$&$\hat{F}^3=$5632& $w^3=7.044399e-01$&$\hat{F}^3=$5623\\
 $w^4=2.086911e-01$&$\hat{F}^4=$6094& $w^4=2.374787e-01$&$\hat{F}^4=$6080\\
 $w^5=8.937330e-19$&$\hat{F}^5=$14800000& $w^5=6.835344e-11$&$\hat{F}^5=$12401053\\
 \hline
\end{tabular}
\end{center}
\end{footnotesize}
\end{table}

\begin{table}[!h]
  \begin{footnotesize}
\begin{center}
\caption{I1. Weights and optimal solution of model LIP-ORA (19) for the cross docking node$-$door assigned costs and cost surplus variables in scenario $\omega$,
 for $\omega\in\Omega_{\mathbb{P}}, \, \mathbb{P}\in\P_{5\%}$}
\label{Tab3.I1.ORA.P5}
  \begin{tabular}{|ccl|ccl|}
    \multicolumn{3}{c}{$\mathbb{P}=1$}&\multicolumn{3}{c}{$\mathbb{P}=2$}\\
  \hline
 $w^1=5.841132e-06 $ &$\hat{F}^1=$5638&$\hat{s}^1=0.0$& $w^1=3.909180e-05$&$\hat{F}^1=$5342&$\hat{s}^1=0.0$\\
  $w^2=6.082577e-07$ &$\hat{F}^2=$5483&$\hat{s}^2=0.0$& $w^2=4.025044e-02$&$\hat{F}^2=$5450&$\hat{s}^2=0.0$\\
 $w^3=7.230249e-01$  &$\hat{F}^3=$5639&$\hat{s}^3=0.0$& $w^3=7.701221e-01$&$\hat{F}^3=$5639&$\hat{s}^3=0.0$\\
  $w^4=2.161435e-01$ &$\hat{F}^4=$5920&$\hat{s}^4=0.0$& $w^4=1.895883e-01$&$\hat{F}^4=$5929&$\hat{s}^4=0.0$\\
  $w^5=2.404819e-12$ &$\hat{F}^5=$6905&$\hat{s}^5=1069$& $w^5=6.835344e-11$&$\hat{F}^5=$7375&$\hat{s}^5=1539$\\
 \hline
 \multicolumn{3}{c}{$\mathbb{P}=3$}&\multicolumn{3}{c}{$\mathbb{P}=4$}\\
 \hline
 $w^1=3.786548e-07$&$\hat{F}^1=$5558&$\hat{s}^1=0.0$& $w^1=7.957016e-06$&$\hat{F}^1=$5187&$\hat{s}^1=0.0$\\
 $w^2=3.658705e-02$&$\hat{F}^2=$5446& $\hat{s}^2=0.0$&$w^2=5.807343e-02$&$\hat{F}^2=$5519& $\hat{s}^2=0.0$\\
 $w^3=7.54214e-01$&$\hat{F}^3=$5632& $\hat{s}^3=0.0$&$w^3=7.044399e-01$&$\hat{F}^3=$5634& $\hat{s}^3=0.0$\\
 $w^4=2.086911e-01$&$\hat{F}^4=$5927& $\hat{s}^4=0.0$&$w^4=2.374787e-01$&$\hat{F}^4=$5920& $\hat{s}^4=0.0$\\
 $w^5=8.937330e-19$&$\hat{F}^5=$7222&$\hat{s}^5=1386$&$w^5=6.835344e-11$&$\hat{F}^5=$6921&$\hat{s}^5=1085$\\
 \hline
   \end{tabular}
 \end{center}
 \end{footnotesize}
\end{table}
\subsection{Ambiguity sets generation for instance I3}\label{sec:dataI3}
The instance is composed of up to 10 origins, 10 destinations, 5 strip doors and 5 stack ones.
In particular, ND is defined by $|\Omega|=10$ equiprobable
scenarios in $|\L|=2$ groups, say $|\Omega^1|=5$ and $|\Omega^2|=5$.

Table \ref{Tab.I3elem} reports the dimensions of the sets in the instance.
The headings are as follows for scenario group $\ell$, for $\ell\in\L$:
$|\Omega^\ell|$, number of scenarios;
$n\M$ ($n\N$), number of origin (destination) nodes;
$n\I$ ($n\J$), number of strip (stack) doors;
$m\H$, cardinality of the set of original uncertain parameters, i.e., $n\M \times n\N$.
\begin{table}[!h]
\begin{center}
 {\small 
\caption{Instance I3. Dimensions of scenario groups in set $\L$}\label{Tab.I3elem}
\begin{tabular}{c|cccccc}
\hline
$\ell$ & $|\Omega^\ell|$ & $n\M$ & $n\N$ & $n\I$ & $n\J$ & $m\H$\\
\hline
1 &  5  &  8  &  8  &   4 &   4  &   64 \\
2 &  5  & 10  & 10  &   5 &   5  &  100 \\
\hline
\end{tabular}
}
\end{center}
\end{table}

Notice that one strip door and one stack door \textit{are} fully capacity reduced for scenario group $\ell=1$,
while the other doors \textit{can} be partially or fully reduced \textit{at random}.
On the other hand, no \textit{given} door \textit{is} fully reduced for scenario group $\ell=2$.
The related non-zero realizations (i.e., commodity volumes $H_{mn}^\omega$) in ND are chosen \textit{at random},
for each scenario $\omega$ in $\Omega^1\cup\Omega^2$.

Similarly to instance $I1$,
Table \ref{Tab2.I3} shows the statistics of $\ell_{\mathbb{P}}^r$ for the candidate members to be included in the ambiguity set $\P$, based on the set of the 80 proximities that are computed in the experiment,
by considering model (B.3).
So, $\P=\P_{7\%}$ corresponds to the first 7-centile $\theta_{7\%}=10.61$ as the \textit{radius};
it is composed of the six consecutively numbered perturbations $\mathbb{P}\in \{1,2,...,6\}$, where the first three are Normal, the fourth and the fifth ones are Gamma and the sixth is Lognormal.
On the other hand, the ambiguity set $\P_{3\%}$ does it for the first 3-centile reduced \textit{radius} $\theta_{3\%}=10.31$;
it is composed of the first three perturbations $\mathbb{P}\in\{1,2,3\}$, so that $\P_{3\%}\subset\P_{7\%}$.
\begin{table}[!h]
\begin{center}
{\small
\caption{Instance I3. Statistics for the candidates for being the ambiguity set}
\label{Tab2.I3}
\begin{tabular}{cccccccc}
  \hline
 Min&1st 3cent&1st 7cent&1st Qu&Median& Mean&3rd Qu&Max\\
  \hline
9.91&10.31&10.61&10.99&11.57&11.68&12.34&14.17\\
  \hline
 \end{tabular}
}
\end{center}
\end{table}
\subsubsection{Instance I3. More details on generating ambiguity sets $\P_{3\%}$ and $\P_{7\%}$ and DRO solutions}\label{sec:I3details}
Table \ref{Tab.I3} summarizes some other statistics for the twenty
proximities $\ell_{\mathbb{P}}^2$;
they are computed under each of the four distribution assumptions in set $\Q$.
The Normal perturbations are the nearest ones to the ND, and the
Lognormal and Gamma ones are also near to the ND.
\begin{table}[!h]
\begin{center}
{\small
\caption{Instance I3. Statistics for \textit{proximity} $\ell_{\mathbb{P}}^r \, \forall \mathbb{P}\in\P'_q, \, q\in\Q$}\label{Tab.I3}
\begin{tabular}{lrrrrrr}
  \hline
\text{Distribution}&  Min & 1st Qu & Median & Mean & 3rd Qu & Max\\
   \hline
  Normal&9.91&10.67&10.98&11.03&11.17&14.17\\
  Weibull&11.02&12.35&12.50&12.68&13.09&14.11\\
  Gamma&10.42&10.99&11.42&11.45&11.77&12.63\\
  Lognormal&10.60&11.01&11.60&11.55&11.75&13.27\\
  \hline
 \end{tabular}
 }
\end{center}
\end{table}

Table \ref{Tab3.I3.RN.P3} shows the weights
and the cross docking node$-$door assignment cost under the scenarios in the ambiguity
set $\P_{3\%}$ as retrieved from the CPLEX incumbent solution for instance I3-RN-3.
The incumbent cost $z_{RN}=8451.89$ is composed of the first stage
infrastructure building cost $C_1= 2082$ and the second stage expected cost $6369.89$ of ambiguity set member $\mathbb{P}=3$,
since $\gamma_{\mathbb{P}}=1$.
 It can be observed the very high cost of the black swan scenarios.
\begin{table}[!h]
\begin{footnotesize}
\begin{center}
\caption{I3. Weights and Cplex incumbent solution of model LIP-RN (8) for the cross docking node$-$door assigned costs under scenario $\omega$, for $\omega\in\Omega_{\mathbb{P}}, \, \mathbb{P}\in\P_{3\%}$}
\label{Tab3.I3.RN.P3}
 \begin{tabular}{|ll|ll|}
    \multicolumn{2}{c}{$\mathbb{P}=1$}&\multicolumn{2}{c}{$\mathbb{P}=2$}\\
  \hline
 $w^1=3.218249e-06$ &$\hat{F}^1=5584$&$w^1=1.228770e-04$&$\hat{F}^1=5775$ \\
  $w^2=3.907031e-02$ &$\hat{F}^2=5623$& $w^2=4.996638e-02$&$\hat{F}^2=5543$\\
 $w^3=2.921724e-01$  &$\hat{F}^3=5624$& $w^3=3.363456e-01 $&$\hat{F}^3=5703$\\
  $w^4=1.687540e-01$ &$\hat{F}^4=5932$& $w^4=1.135651e-01$&$\hat{F}^4=6014$ \\
  $w^5=1.684348e-10$ &$\hat{F}^5=9803406$&$w^5= 3.508562e-12 $&$\hat{F}^5=14200432$\\
 $w^6=3.727733e-08 $ &$\hat{F}^6=9405025$& $w^6= 2.333095e-09$&$\hat{F}^6=9203768$ \\
  $w^7=9.015966e-03$ &$\hat{F}^7=6940$& $w^7=3.080624e-03$&$\hat{F}^7=6888$ \\
 $w^8= 4.071194e-01$  &$\hat{F}^8=7111$& $w^8= 4.551538e-01$&$\hat{F}^8=7159$\\
 $w^9= 8.386457e-02 $ &$\hat{F}^9=7703$& $w^9=4.176560e-02$&$\hat{F}^9=7685$ \\
 $w^{10}= 9.389927e-21$ &$\hat{F}^{10}=15003648$&$w^{10}=4.330379e-18$&$\hat{F}^{10}=15001678$\\
   \hline
    \multicolumn{2}{c}{$\mathbb{P}=3$}&\\
     \cline{1-2}
   $w^1= 1.020567e-07$&$\hat{F}^1=13200432$\\
   $w^2= 7.065075e-02$&$\hat{F}^2=5501$\\
   $w^3=2.021934e-01$&$\hat{F}^3=5592$\\
   $w^4=2.271557e-01$&$\hat{F}^4=5948$\\
   $w^5=3.077295e-18$&$\hat{F}^5=9801489$\\
   $w^6=5.518698e-11$&$\hat{F}^6=18801824$\\
   $w^7= 1.757278e-03$&$\hat{F}^7=6637$\\
   $w^8=4.599281e-01$&$\hat{F}^8=7201$\\
   $w^9= 3.831462e-02$&$\hat{F}^9=7539$ \\
$w^{10}=2.904815e-27$&$\hat{F}^{10}=15001564$\\
   \cline{1-2}
\end{tabular}
\end{center}
\end{footnotesize}
\end{table}

Table \ref{I3-ORA-3.CPX-sol} shows some results retrieved from the incumbent solution of models LIP-RN (8) and LIP-ORA (19),
as the solution cost $\hat{F}^\omega$ of the cross docking node$-$door assignment, the total
cost $\hat{C}_{12,\mathbb{P}=3}^\omega$ and the cost surplus
variables $s^\omega$. For the ORA parameters' triplex ($\iota=9600$, $\overline{s}=2500$, $\rho=3000$),
the incumbent cost $z_{ORA}=8526.87$ is composed of the first stage infrastructure building cost $\hat{C}_1=2082$ and the second stage  expected cost 6444.87 in ambiguity set member $\mathbb{P}=3$, since $\gamma_\mathbb{P}=1$, for $\mathbb{P}\in\P_{3\%}$.
The total cost $\hat{C}_{12,\mathbb{P}=3}^\omega$, see Section 4.1 of the main body, is composed of $\hat{C}_1$ plus the second stage cost $\hat{F}^\omega$ for scenario $\omega$.
It can be seen that the $s$-variables in the incumbent RN solution have zero-values for all scenarios in the ORA version, but scenario $\omega=10$,
which value is smaller than bound $\overline{s}$ in ORA constraint (18b).
\begin{table}[!h]
\begin{center}
    {\scriptsize 
\caption{I3. Cplex incumbent solution of models LIP-RN (8) and LIP-ORA (19), in particular cost surplus variables for the selected member $\mathbb{P}=3:~\gamma_{3}=1$ in ambiguity set $\P_{3\%}$}
\label{I3-ORA-3.CPX-sol}
\begin{tabular}{l|ccc|ccc}
  \hline
  & \multicolumn{3}{c}{I3-RN-3}&\multicolumn{3}{|c}{I3-ORA-3}\\
   \hline
  &  $\hat{F}^\omega$ & $C_{12,\mathbb{P}=3}^\omega$&$\hat{s}^\omega$&$\hat{F}^\omega$ & $C_{12,\mathbb{P}=3}^\omega$&$\hat{s}^\omega$\\
\hline
$w^1= 1.020567e-07$&13200432&13202514&13192914&5603&7685&0\\
$w^2= 7.065075e-02$&5501&7583&0&5466&7548&0\\
$w^3=2.021934e-01$&5592&7674&0&5602&7684&0\\
$w^4=2.271557e-01$&5948&8030&0&5910&7992&0\\
$w^5=3.077295e-18$&9801489&9803571&9793971&7480&9562&0\\
$w^6=5.518698e-11$&18801824&18803906&18794306&6673&8755&0\\
$w^7= 1.757278e-03$&6637&8719&0&6658&8740&0\\
$w^8=4.599281e-01$&7201&9283&0&7142&9224&0\\
$w^9= 3.831462e-02$&7539&9621&21&7491&9573&0 \\
$w^{10}=2.904815e-27$&15001564&1553646&14994046&9405&11487&1887\\
  \hline
\end{tabular}
}
\end{center}
\end{table}

Table \ref{I3-3.No-Violation} shows the matheuristic (see Sections 3.3 and 4.3) incumbent solution of models LIP-RN (8) and LIP-ORA (19).
The solution cost ($\overline{z}_H= 8320.96, \hat{C}_{1}=1849$) for I3-RN-3 violates the non-imposed upper bound $\overline{s}=2500$ in the black swan scenarios; it does not happen for I3-ORA-3.
\begin{table}[!h]
\begin{center}
    {\scriptsize 
\caption{I3. Matheuristic incumbent solution of models LIP-RN (8) and LIP-ORA (19), in particular, the cost surplus variables for the selected member $\mathbb{P}=3:~\gamma_{3}=1$ in ambiguity set $\P_{3\%}$}\label{I3-3.No-Violation}
\begin{tabular}{l|ccc|ccc}
  \hline
  & \multicolumn{3}{c}{I3-RN-3}&\multicolumn{3}{|c}{I3-ORA-3}\\
   \hline
  &  $\hat{F}_{H}^\omega$ & $\hat{C}_{12,\mathbb{P}=3}^\omega$&$\hat{s}^\omega$&$\hat{F}_{H}^\omega$ & $\hat{C}_{12,\mathbb{P}=3}^\omega$&$\hat{s}^\omega$\\
\hline
$w^1= 1.020567e-07$&14800000&14801849&14792276&5640&7722&0\\
  $w^2= 7.065075e-02$&5490&7339&0&5466&7548&0\\
$w^3=2.021934e-01$&5617&7466&0&5592&7674&0\\
$w^4=2.271557e-01$&5998&7847&0&5910&7792&0\\
$w^5=3.077295e-18$&14800000&14801849&14792276&7601&9683&110\\
$w^6=5.518698e-11$&22000000&22001849&21992276&7027&9109&0\\
$w^7= 1.757278e-03$&6738&8587&0&6648&8730&0\\
$w^8=4.599281e-01$&7159&9008&0&7099&9181&0\\
$w^9= 3.831462e-02$&7549&9398&0&7459&9541&0 \\
$w^{10}=2.904815e-27$&22000000&22001849&21992276&9344&11426&1853\\
  \hline
\end{tabular}
}
\end{center}
\end{table}

In a similar way to $\P_{3\%}$, Table \ref{I3-6.No-Violation} shows the matheuristic incumbent solution of models LIP-RN (8) and LIP-ORA (19) in ambiguity set $\P_{7\%}$.
The solution ($\overline{z}_H= 8327.80$, $\hat{C}_{1}=1849$) for I3-RN-6 violates the non-imposed upper bound $\overline{s}=2500$;
it does not happen for I3-ORA-6.
\begin{table}[!h]
\begin{center}
    {\scriptsize 
\caption{I3. Matheuristic incumbent solution of models LIP-RN (8) and LIP-ORA (19), in particular, the cost surplus variables for the selected member $\mathbb{P}=3:~\gamma_{3}=1$ in ambiguity set $\P_{7\%}$}\label{I3-6.No-Violation}
\begin{tabular}{l|ccc|ccc}
  \hline
  & \multicolumn{3}{c}{I3-RN-6}&\multicolumn{3}{|c}{I3-ORA-6}\\
   \hline
  &  $\hat{F}_{H}^\omega$ & $\hat{C}_{12,\mathbb{P}=3}^\omega$&$\hat{s}^\omega$&$\hat{F}_{H}^\omega$ & $\hat{C}_{12,\mathbb{P}=3}^\omega$&$\hat{s}^\omega$\\
\hline
$w^1= 1.020567e-07$&14800000&14801849&14792276&5334&7623&0\\
  $w^2= 7.065075e-02$&5503&7532&0&5466&7555&0\\
$w^3=2.021934e-01$&5634&7483&0&5592&7681&0\\
$w^4=2.271557e-01$&5998&7844&0&5910&7999&0\\
$w^5=3.077295e-18$&14800000&14801849&14792276&7629&9718&170\\
$w^6=5.518698e-11$&22000000&22001849&21992276&6822&8911&0\\
$w^7= 1.757278e-03$&6689&8538&0&6710&8799&0\\
$w^8=4.599281e-01$&7159&9008&0&7099&9188&0\\
$w^9= 3.831462e-02$&7595&9444&0&7459&9548&0 \\
$w^{10}=2.904815e-27$&22000000&22001849&21992276&9224&11313&1765\\
  \hline
\end{tabular}
}
\end{center}
\end{table}
\subsection{Ambiguity sets generation for instance I7}\label{sec:dataI7}
The instance is composed of up to 20 origins, 20 destinations, 10 strip doors and 10 stack ones.
In particular, ND is defined by $|\Omega|=20$ equiprobable scenarios in $|\L|=4$ groups.
Table \ref{Tab.I7elem} reports the dimensions of the sets in the instance.
Notice that 6 strip doors and 6 stack doors are fully capacity reduced in the scenarios group $\ell=1$, 5 strip doors and 5 stack doors in the case of scenario group $\ell=2$, and 4 strip doors and 4 stack doors for scenario group $\ell=3$. The other doors can be partially or fully reduced \textit{at random}.
\begin{table}[!h]
\begin{center}
 {\small
\caption{Instance I7. Dimensions of scenario groups in set $\L$}\label{Tab.I7elem}
\begin{tabular}{c|cccccc}
\hline
$\ell$ & $|\Omega^\ell|$ & $n\M$ & $n\N$ & $n\I$ & $n\J$ & $m\H$\\
\hline
1 &  5  &  8  &  8  &   4 &   4  &   64\\
2 &  5  & 10  & 10  &   5 &   5  &  100\\
3 &  5  & 15  & 15  &   6 &   6  &  225\\
4 &  5  & 20  & 20  &  10 &  10  &  400\\
\hline
\end{tabular}
}
\end{center}
\end{table}

Following the same scheme, the
candidate ambiguity set is built with $|\P'_q|=20$ and $|\Q|=4$, and
Table \ref{Tab2.I7} shows the statistics of $\ell_{\mathbb{P}}^r$ for the candidate members to be included in the ambiguity set $\P$,
based on the set of the 80 proximities and by considering model (B.3).
So, $\P=\P_{7\%}$ corresponds to the first 7-centile $\theta_{7\%}=12.58$ as the \textit{radius};
it is composed of the six consecutively numbered perturbations $\mathbb{P}\in \{1,2,...,6\}$,
where the second and the third ones are Lognormal, the sixth is Weibull  and, the rest are Gamma.
On the other hand, the ambiguity set $\P_{3\%}$ does it for the first 3-centile reduced \textit{radius} $\theta_{3\%}=12.29$;
it is composed of the first three perturbations $\mathbb{P}\in\{1,2,3\}$, so that $\P_{3\%}\subset\P_{7\%}$,
of which the first one is Gamma and the other two Lognormal.
\begin{table}[h!]
\begin{center}
{\small
\caption{Instance I7. Statistics for the candidate for being the ambiguity sets}
\label{Tab2.I7}
\begin{tabular}{cccccccc}
  \hline
 Min&1st 3cent&1st 7cent&1st Qu&Median& Mean&3rd Qu&Max\\
  \hline
11.80&12.29&12.58&13.30&13.85&13.93&14.40&17.53\\
  \hline
 \end{tabular}
 }
\end{center}
\end{table}
\subsubsection{Instance I7. More details on generating ambiguity sets $\P_{3\%}$ and $\P_{7\%}$ and DRO solutions}\label{sec:I7details}
Table \ref{Tab.I7} summarizes some other statistics for the twenty
proximities $\ell_{\mathbb{P}}^2$ that are computed under each of the four probability distribution assumptions in set $\Q$.
As can be seen the Gamma and Lognormal  perturbations are the nearest to the ND.
\begin{table}[!h]
\begin{center}
{\small
\caption{I7. Statistics for \textit{proximity} $\ell_{\mathbb{P}}^r \, \forall \mathbb{P}\in\P'_q, \, q\in\Q$}\label{Tab.I7}
\begin{tabular}{lrrrrrr}
  \hline
Distribution&  Min & 1st Qu & Median & Mean & 3rd Qu & Max\\
   \hline
  Normal&13.87&14.16&14.70&14.98&15.12&17.53\\
  Weibull&12.53&13.52&13.71&13.87&14.12&15.88\\
  Gamma&11.80&12.80&13.34&13.40&14.01&14.88\\
  Lognormal&12.19&13.20&13.55&13.48&13.80&14.41\\
  \hline
 \end{tabular}
 }
\end{center}
\end{table}

Table \ref{I7-3.No-Violation} shows the results provided by the matheuristic for the instances LIP-RN-3 (8) and LIP-ORA-3 (19),
where the ORA parameters' triplex ($\iota=38800$, $\overline{s}=7900$, $\rho=5000$) is considered in ambiguity set $\P_{3\%}$.
It can be observed that the incumbent solution cost for I7-RN-3 in member $\mathbb{P}=3$
($\overline{z}_H=22402.18, \hat{C}_{1}=7157$) does not satisfy the non-imposed upper bound $\overline{s}=7900$.
It can be observed that the black swan scenarios have not any more so negative
implications in the ORA-solution ($\overline{z}_H=22402.69, \hat{C}_{1}=7284$) with a minimal increase in $z$-cost.
\begin{table}[!h]
    {\scriptsize 
\caption{I7. Matheuristic incumbent solution of instances LIP-RN-3 (8) and LIP-ORA-3 (19), in particular, the cost surplus variables for the
  selected members $\mathbb{P}=3:~\gamma_{3}=1$ (RN) and $\mathbb{P}=1:~\gamma_{1}=1$ (ORA) for ambiguity set $\P_{3\%}$}\label{I7-3.No-Violation}
\begin{tabular}{l|ccc|l|ccc}
  \hline
  \multicolumn{4}{c}{I7-RN-3}&\multicolumn{4}{|c}{I7-ORA-3}\\
   \hline
  weight for $\mathbb{P}=3$ &  $\hat{F}_{H}^\omega$ & $\hat{C}_{12,\mathbb{P}=3}^\omega$&$\hat{s}^\omega$&weights $\mathbb{P}=1$&$\hat{F}_{H}^\omega$ & $\hat{C}_{12,\mathbb{P}=1}^\omega$&$\hat{s}^\omega$\\
\hline
$w^1= 5.890609e-05$&5026&12183&0&5.034938e-05&5120&12404&0\\
  $w^2= 1.113443e-01 $&5372&12527&0&1.522775e-01 &5237&12521&0\\
$w^3= 1.140110e-01$&5519&12676&0& 6.236860e-02 &5663&12947&0\\
$w^4=2.458577e-02$&5791&12948&0& 3.530357e-02&5819&13103&0\\
$w^5=1.597404e-12 $&6505&13662&0&   9.910015e-13&6496&13780&0\\
$w^6= 2.151672e-13$&5808&12965&0&  1.085449e-11&5905&13189&0\\
$w^7= 3.399928e-02 $&6425&13582&0& 1.473375e-02 &6323&13607&0\\
$w^8=  2.156317e-01$&6612&13769&0& 2.348509e-01 &6628&13912&0\\
$w^9=  3.690257e-04 $&7028&14185&0& 4.153145e-04 &6848&14132&0 \\
$w^{10}=  1.182285e-18$&7727&14884&0&  1.341637e-28&8259&15543&0\\
$w^{11}= 1.717728e-17 $&14252&21409&0&2.357762e-16 &14220&21504&0\\
$w^{12}=  1.142020e-03 $&14945&22102&0& 5.588144e-04 &15053&22337&0\\
$w^{13}= 2.488110e-01$&16000&23157&0& 2.492375e-01&15863&23147&0\\
$w^{14}= 4.698466e-05$&16910&24067&0& 2.036854e-04 &16706&23990&0\\
$w^{15}= 1.111944e-47$&18697&25854&0&  1.124869e-51 &18818&26102&0\\
$w^{16}= 1.445404e-36$&29396&36553&0&  5.365765e-38 &29135&36419&0\\
$w^{17}= 2.012968e-09$&31040&38197&0&  1.650646e-09&30573&37857&0\\
$w^{18}=   2.500000e-01$&32918&40705&1275& 2.500000e-01  &32557&39861&1061
\\
$w^{19}= 3.104054e-08$&22000000&22007157& 21968357&1.794206e-09 &34757&42041&3241
\\
$w^{20}=  7.772709e-97$&22000000&22007157& 21968357& 1.737634e-72&39408&46692&7892
\\
  \hline
\end{tabular}
}
\end{table}

For the same ORA parameter's triplex, Table \ref{I7-6.No-Violation} shows the matheuristic incumbent solution of instances LIP-RN-6 (8) and LIP-ORA-6 (19).
In particular, it shows the cost surplus variables for the selected member $\mathbb{P}=6$, since $\gamma_{\mathbb{P}}=1$, in ambiguity set $\P_{7\%}$,
for the instances I7-RN-6 ($\overline{z}_H= 23700.67$, $\hat{C}_{1}=7784$) and I7-ORA-6 ($\overline{z}_H= 23939.31$, $\hat{C}_{1}=8023$).
It can be observed that the black swan scenarios in model LIP-RN (8) have not any more so negative implications in the ORA-solution with a minimal increase in cost.
\begin{table}[h!]
    {\scriptsize 
\caption{I7. Matheuristic incumbent solutions of instances LIP-RN-6 (8) and
  LIP-ORA-6 (19), in particular, the cost surplus variables for the
  selected member $\mathbb{P}=6:~\gamma_{6}=1$ in ambiguity set $\P_{7\%}$}\label{I7-6.No-Violation}
\begin{tabular}{l|ccc|ccc}
  \hline
  &\multicolumn{3}{c}{I7-RN-6}&\multicolumn{3}{|c}{I7-ORA-6}\\
   \hline
  weight for $\mathbb{P}=6$ &  $\hat{F}_{H}^\omega$ & $\hat{C}_{12,\mathbb{P}=3}^\omega$&$\hat{s}^\omega$&$\hat{F}_{H}^\omega$ & $\hat{C}_{12,\mathbb{P}=1}^\omega$&$\hat{s}^\omega$\\
\hline
$w^1= 6.231762e-05$&5185&12969&0&5185&13208&0\\
  $w^2= 5.578641e-02$&5350&13134&0&5350&13373&0\\
$w^3=  1.190980e-01$&5551&13335&0&5551&13574&0\\
$w^4=7.505325e-02$&5827&13611&0&5839&13862&0\\
$w^5=1.286116e-15 $&6703&14487&0&6703&14726&0\\
$w^6=  2.679631e-08 $&5954&13738&0&5954&13977&0\\
$w^7= 1.521364e-03$&6336&14120&0&6336&14359&0\\
$w^8= 4.956362e-02 $&6545&14329&0&6545&14568&0\\
$w^9=  1.989150e-01 $&7015&14799&0&7015&15038&0\\
$w^{10}=  1.380043e-31 $&7774&15558&0&7811&15834&0\\
$w^{11}=  8.606838e-18 $&14295&22079&0&14298&22321&0\\
$w^{12}=  1.413645e-07 $&15190&22974&0&15194&23217&0\\
$w^{13}= 1.715785e-03  $&15914&23698&0&15900&23923&0\\
$w^{14}=  2.482841e-01$&16608&24392&0 &16597&24620&0\\
$w^{15}= 3.971946e-47$&18831&26615&0 &18889&26912&0\\
$w^{16}=  1.731425e-38$&29447& 37231
&0& 29289
& 37312
&0\\
$w^{17}=   1.524912e-16$& 31342& 39126
&326
& 31142
&39165
&365
\\
$w^{18}=   1.520760e-05 $& 33018
& 40802
&2002
 & 32632
&40655
&1855
\\
$w^{19}= 2.499848e-01 $& 34557
& 42341
& 3541
 & 34567
& 42590
&3790
\\
$w^{20}=  8.559141e-75$&22000000&22007784
& 21968984
&38521
& 46544
& 7744
\\
  \hline
\end{tabular}
}
\end{table}

\end{document}